\documentclass[11pt]{article}

\newcommand{\keywords}[1]{\par\medskip\noindent\textbf{Keywords:} #1\par}

\usepackage{amsmath,amsthm,amsopn,amsfonts,amssymb}
\usepackage{mathrsfs,a4wide}
\usepackage{indentfirst}
\usepackage{xcolor}

\usepackage{graphicx,graphics,color}
\newtheorem{thm}{Theorem}[section]
\newtheorem{cor}[thm]{Corollary}
\newtheorem{Lemma}[thm]{Lemma}
\newtheorem{prop}[thm]{Proposition}

\newtheorem{ppt}[thm]{Properties}

\theoremstyle{definition}
\newtheorem{defn}[thm]{Definition}
\theoremstyle{remark}
\newtheorem{rem}[thm]{Remark}

\numberwithin{equation}{section}

\newcommand{\norm}[1]{\left\Vert#1\right\Vert}

\newcommand{\N}{\mathbb{N}}
\newcommand{\R}{\mathbb{R}}
\newcommand{\Z}{\mathbb{Z}}

\begin{document}

\tolerance=1000
\emergencystretch=2em

\title{The first non-trivial Fu\v{c}\'\i k curve for the sublinear Dirichlet Laplacian}
\author{Elie Alhelou \and Gisella Croce \and Giovanni Pisante}
\date{}

\maketitle

\begin{abstract}
Let $\Omega\subset\R^N$ be a bounded open connected set, and let $q \in (1,2)$. We consider the sublinear $q$-eigenvalue problem for the Dirichlet Laplacian
$$
    -\Delta u = \lambda \norm{u}_{L^q(\Omega)}^{2-q} |u|^{q-2}u, \qquad u\in H^1_0(\Omega)\setminus\{0\},
$$
whose first eigenvalue $\lambda_1(\Omega)$ is simple, with eigenfunctions of constant sign. The definition of a second eigenvalue is {not unique}. Following the approach introduced by Molle and Passaseo for the linear Laplacian, we construct the first non-trivial curve of the corresponding Fu\v{c}\'\i k spectrum, and show that it has a unique diagonal point $(\gamma(\Omega),\gamma(\Omega))$. The Fu\v{c}\'\i k-selected second value $\lambda_{2}^F(\Omega) = 2^{1 - \frac{2}{q}}\gamma(\Omega)$ is strictly larger than $\lambda_1(\Omega)$, admits a two-phase variational characterisation, and satisfies a sharp Hong--Krahn--Szeg\H{o} inequality: among sets of given measure it is minimised by the union of two equal balls. The selected solution at the diagonal point is a classical $(2,q)$-eigenfunction exactly when its nodal $L^q$-masses are balanced. We finally show that on balls in dimensions $N\geq2$ every Fu\v{c}\'\i k-selected solution is foliated Schwarz symmetric.
\end{abstract}


\tableofcontents

\section{Introduction}

Given a bounded open connected set $\Omega\subset\R^N$ and $q\in\,(1,2)$, we study the non-linear $q$-eigenvalue problem for the Dirichlet Laplacian
\begin{equation*}
    \left\{
    \begin{array}{cl}
        -\Delta u = \lambda\norm{u}_{L^q(\Omega)}^{2-q}|u|^{q-2}u, & \textnormal{ in } \Omega \\
        u = 0, & \textnormal{ on } \partial\Omega
    \end{array}
    \right.
\end{equation*}
whose solutions $(\lambda,u)\in\R^*_+\times H_0^1(\Omega)\setminus\{0\}$ are the critical points of the Rayleigh-type quotient
$$
    u \; \longmapsto \; \frac{\displaystyle\int_\Omega |\nabla u|^2 \,dx}{\left(\displaystyle \int_\Omega |u|^q \,dx\right)^{\frac{2}{q}}}\,.
$$
Problems of this kind, often referred to as {sublinear Lane--Emden} or $(p,q)$-eigenvalue problems, have attracted considerable attention in recent years: we refer to Lindqvist \cite{Lindqvist} for a general overview of non-linear eigenvalue problems, and to Franzina and Lamberti \cite{Franzina_Lamberti} for the $(p,q)$-eigenvalue problem.

The first eigenvalue $\lambda_1(\Omega)$, defined by minimising the quotient above, enjoys the familiar properties of the case $q=2$: it is simple, its eigenfunctions have constant sign, and a Faber--Krahn inequality holds.
On the other hand, even the mere definition of a “second eigenvalue” is  a delicate question: indeed, the quantity
$$
    \lambda_2^\sigma(\Omega) = \inf \{\lambda\in\sigma \,;\, \lambda > \lambda_1(\Omega)\}
$$
is defined as an infimum whenever the set is non-empty. It is strictly larger than $\lambda_1(\Omega)$, and is itself an eigenvalue, if $\lambda_1(\Omega)$ is isolated in {the spectrum} $\sigma$, as established by Brasco, De Philippis and Franzina \cite[Theorem B and Theorem 6.3]{Brasco_DePhilippis_Franzina} for $N\geq2$ on bounded open sets of class $C^1$ with finitely many connected components, and on Lipschitz open sets for $q$ in a range $(q_\Omega,2)$ depending on the Lipschitz constant.

It is therefore useful to distinguish several natural candidates for a second level, independently of the above definition of $\lambda_2^\sigma(\Omega)$, given by different variational characterisations. We set
$$
    \mathcal{M}_q := \{u \in H_0^1(\Omega) : \norm{u}_{L^q(\Omega)} = 1\}, \qquad E(u) := \int_\Omega |\nabla u|^2 \,dx,
$$
and let $e_1\in\mathcal{M}_q$ be the positive first $(2,q)$-eigenfunction. The second
Ljusternik--Schnirelman value is defined as
$$
    \lambda_2^{LS}(\Omega) := \inf_{A\in\mathscr{A}_2}\ \sup_{u\in A}\ E(u),
$$
where $\mathscr{A}_2$ is the family of compact symmetric subsets of $\mathcal{M}_q$ having Krasnosel'skii genus at least two {(the genus of a closed symmetric set $A\subset\mathcal{M}_q$ being the smallest $k\in\N$ such that there exists an odd continuous map $A\to\R^k\setminus\{0\}$, and $+\infty$ if no such $k$ exists)}, as has been stated in \cite{Franzina_Lamberti}. We can also adapt the mountain-pass value between the two first $(2,q)$-eigenfunctions, as seen in \cite{Cuesta_deFigueiredo_Gossez} for the $p$-Laplacian eigenvalue problem, as follows
$$
    \lambda_2^{MP}(\Omega) := \inf_{\gamma\in\mathscr{P}} \ \max_{t\in[-1,1]} \ E(\gamma(t)), \qquad \mathscr{P} := \{\gamma\in C([-1,1],\mathcal M_q) \,;\, \gamma(-1) = -e_1,\ \gamma(1) = e_1\}.
$$
We may also consider the least sign-changing critical level
$$
    \lambda_2^{nod}(\Omega) := \inf\{E(u):u\in\mathcal{M}_q,\ E'|_{\mathcal{M}_q}(u) = 0, \ u^+\not\equiv0,\ u^-\not\equiv0\},
$$
whenever this set is non-empty, where $u^+$ (respectively $u^-$) is the positive (respectively negative) part of $u$.  In the linear case $q=2$, all these quantities coincide with the usual second eigenvalue $\lambda_2^\sigma(\Omega)$, that is in fact well-defined. For $1<q<2$, no such general identification is known.

In this paper, we propose a definition inspired by the work of Molle and Passaseo \cite{Molle_Passaseo} for the Laplacian. The Fu\v{c}\'\i k spectrum $\Sigma$ associated to our problem is the set of pairs $(\alpha,\beta)\in\R^2$ for which the equation
$$
    -\Delta u = -\alpha \norm{u^-}_{L^q(\Omega)}^{2-q}(u^-)^{q-1}+\beta \norm{u^+}_{L^q(\Omega)}^{2-q}(u^+)^{q-1}
$$
admits a non-trivial solution $u\in H^1_0(\Omega)$. It contains the trivial lines $\{\lambda_1(\Omega)\}\times\R$ and $\R\times\{\lambda_1(\Omega)\}$, and every eigenvalue $\lambda \in \sigma\setminus\{\lambda_1(\Omega)\}$ generates points of $\Sigma$ off the trivial lines. Following \cite{Molle_Passaseo}, for every $\beta > \lambda_1(\Omega),$ we construct, by minimisation of a family of regularised functionals $f_{\beta,\varepsilon}$ and a limiting procedure as $\varepsilon\to0$, a value
$$
    \alpha_\beta = \min\left\{\int_\Omega |\nabla u^-|^2\,dx \,;\, u\in H_0^1(\Omega), \norm{u^+}_{L^q(\Omega)} = \norm{u^-}_{L^q(\Omega)} = 1 , \int_\Omega |\nabla u^+|^2\,dx=\beta \right\},
$$
such that $(\alpha_\beta,\beta)\in\Sigma$. The curve $\beta \mapsto \alpha_\beta$ is the \textit{first non-trivial curve} of $\Sigma$: it is minimal among the points of $\Sigma$ with $\alpha,\beta>\lambda_1(\Omega)$, continuous and strictly decreasing, it tends to $\lambda_1(\Omega)$ as $\beta\to+\infty$, and blows up as $\beta\to\lambda_1(\Omega)^+$.
The general strategy follows \cite{Molle_Passaseo}, but the non-local coupling induced by the global $L^q$ norm and the separate normalisation of the two nodal phases require different arguments, in particular in the support-modification steps, where the sublinearity $q<2$ enters: adding a small disjoint bump of size $s$ to a normalised function costs $O(s^2)$ in energy and gains $O(s^q)$ in the constraint (see the proof of Lemma \ref{ModificationOfSupports}).
Moreover, the behaviour of the curve as $\beta\to\lambda_1(\Omega)^+$ is obtained without using the isolation of $\lambda_1(\Omega)$ in $\Sigma$, for the whole range $1<q<2$.

By the strict monotonicity and the asymptotics of the curve, its intersection with the diagonal $\{\alpha=\beta\}$ consists of a single point $(\gamma(\Omega),\gamma(\Omega))\in \Sigma$ such that $\alpha_{\gamma(\Omega)} = \gamma(\Omega)$. We show that a Fu\v{c}\'\i k-selected solution $u_{\gamma}$ at this point becomes a solution of \eqref{Eigen2q} with $\lambda=2^{1-\frac2q}\gamma(\Omega)$ if and only if the following additional balance condition $\norm{u_\gamma^+}_{L^q(\Omega)} = \norm{u_\gamma^-}_{L^q(\Omega)}$ holds.

Defining $\lambda_{2}^F(\Omega) := 2^{1 - \frac{2}{q}} \gamma(\Omega)$, we prove a Hong--Krahn--Szeg\H{o} type inequality. We first show that $\lambda_2^F(\Omega)>\lambda_1(\Omega)$ for every bounded connected open set $\Omega$ (Proposition \ref{lambda1_less_lambda2F}), and that $\gamma(\Omega)$ admits a two-phase characterisation: it is the infimum of $\max\{E(v),E(w)\}$, where $E(u)=\int_\Omega|\nabla u|^2\,dx$,
over the pairs of non-negative functions $v,w\in\mathcal M_q$ with disjoint supports (Proposition \ref{Gamma = gamma}).
At fixed volume $V$, the infimum of $\lambda_{2}^F$ among connected open sets equals $2^{1 - \frac{2}{q}} \lambda_1(B_{V/2})$ and is approached by dumbbell domains. After extending the two-phase definition to possibly disconnected sets, we show that the minimum is attained by the union of two equal balls of measure $\frac{V}{2}$. More precisely, for a bounded open set $D$ we set $\lambda_2^F(D):=2^{1-\frac2q}\Gamma(D)$, where $\Gamma(D)$ is the two-phase infimum above; then $\lambda_2^F(D)\geq2^{1-\frac2q}\lambda_1(B_{V/2})$ whenever $|D|=V$, with equality when $D$ is the union of two disjoint balls of measure $\frac V2$. Conversely, equality implies that $D$ agrees with such a union up to a set of measure zero (Corollary \ref{GeneralShapeOptimisationGamma}).

We also study a symmetry problem. When $\Omega=B\subset\R^N$, $N\geq2$, is a ball, every Fu\v{c}\'\i k-selected solution obtained by the regularised minimisation is \textit{foliated Schwarz symmetric}. This means that the solution depends only on $|x|$ and on the angle $\theta=\arccos\big(\frac{x}{|x|}\cdot e\big)\in[0,\pi]$ from a fixed direction $e\in\mathbb S^{N-1}$, and is non-increasing in $\theta$. The proof reproduces the two-point argument of Gir\~ao and Weth \cite{Girao_Weth}.
The nodal balance problem for diagonal Fu\v{c}\'\i k-selected solutions, and their identification with classical $(2,q)$-eigenfunctions, will be studied for $q$ close to $2$ in the companion paper \cite{PartII}, which is in preparation.

The paper is organised as follows. Section~2 collects the basic properties of the first $(2,q)$-eigenvalue and the regularity of the $(2,q)$-eigenfunctions. Section~3 contains the construction of the first non-trivial Fu\v{c}\'\i k curve, its qualitative properties and the regularity of the Fu\v{c}\'\i k $(2,q)$-eigenfunctions. Section~4 establishes the two-phase characterisation of $\gamma(\Omega)$ and the Hong--Krahn--Szeg\H{o} inequality, and Section~5 the foliated Schwarz symmetry of the Fu\v{c}\'\i k-selected solutions on the ball.

\section{The $(2,q)$-eigenvalue problem}

Let $\Omega$ be a bounded open set in $\R^N$, and let $q\in(1,2^*)$, where $2^*$ is the Sobolev conjugate of $2$ given by
$$
    2^* = \left\{\begin{array}{cl}
        \frac{2N}{N-2}, & \textnormal{ if } N > 2 \\
        +\infty, & \textnormal{ if } N \leq 2.
    \end{array}\right.
$$
We consider the $q$-eigenvalue problem for the Laplacian with Dirichlet boundary conditions, consisting of finding $\lambda\in\R_+^*$ and $u\in H_0^{1}(\Omega) \setminus \{0\}$ such that
\begin{equation}\label{Eigen2q}
    \left\{
    \begin{array}{cl}
        -\Delta u = \lambda\norm{u}_{L^q(\Omega)}^{2-q}|u|^{q-2}u & \textnormal{ in } \Omega \\
        u = 0 & \textnormal{ on } \partial\Omega.
    \end{array}
    \right.
\end{equation}
For simplicity, we will refer to a $q$-eigenvalue of the Laplacian as a $(2,q)$-eigenvalue. We denote by $\sigma$ the \textbf{$(2,q)$-spectrum}, that is the set of all $(2,q)$-eigenvalues $\lambda$. The set $\sigma$ contains a smallest element, denoted by $\lambda_1(\Omega)$ and called the \textbf{first $(2,q)$-eigenvalue}. It admits the following variational characterisation
\begin{equation}\label{FirstEV}
    \lambda_1(\Omega) = \inf\left\{\frac{\norm{\nabla u}^2_{[L^2(\Omega)]^N}}{\norm{u}^2_{L^q(\Omega)}} ; u\in H_0^{1}(\Omega)\setminus\{0\}\right\}.
\end{equation}
The infimum is achieved by a function $u\in H_0^{1}(\Omega)$, called a \textbf{first $(2,q)$-eigenfunction}.

\subsection{Properties of the first $(2,q)$-eigenvalue}

In the following, we will always assume that $\Omega$ is a connected open set of finite measure. We recall a handful of interesting results from \cite{Franzina_Lamberti}.

\begin{prop}\label{FirstEigenFunctionsArePositive}
    Let $q\in (1,2^*)$. If $u\geq 0$ is a $(2,q)$-eigenfunction associated with a $(2,q)$-eigenvalue $\lambda\in\sigma$, then $u>0$ a.e. in $\Omega$. Hence, every first $(2,q)$-eigenfunction associated with $\lambda_1(\Omega)$ must be of constant sign a.e. on $\Omega$.
\end{prop}
For the proof see \cite[Theorem 2.2 and Corollary 2.3]{Franzina_Lamberti}.

\begin{prop}\label{Simplicity_lambda1}
    Let $q\in (1,2]$. Then, the first $(2,q)$-eigenvalue $\lambda_1(\Omega)$ is simple: the associated eigenspace is one-dimensional.
\end{prop}
For the proof see \cite[Theorem 3.1]{Franzina_Lamberti}.

\begin{prop}\label{PositiveEigenfunctionsAreFirst}
    Let $q\in (1,2]$. If $u$ is a $(2,q)$-eigenfunction of constant sign a.e. in $\Omega$, then $u$ is associated to $\lambda_1(\Omega)$.
\end{prop}
For the proof see \cite[Theorem 4.1]{Franzina_Lamberti},  stated there for positive eigenfunctions (if $u\leq0$, apply it to $-u$).

\begin{cor}\label{StrictMonotonicityEV1}
    Let $A$ and $B$ be two subsets of $\R^N$ such that $B$ is connected. If $A\subset B$, then $\lambda_1(A)\geq\lambda_1(B)$. Moreover, if $|B\setminus A| \neq 0$, then $\lambda_1(A)>\lambda_1(B)$.
\end{cor}

We may also show a scaling Lemma for the first $(2,q)$-eigenfunction, as follows.

\begin{Lemma}\label{Homogeneity}
    Let $q \in (1,2^*)$. The first eigenvalue $\lambda_1(\Omega)$ satisfies the following  homogeneity property:
    $$
        \forall t>0 \, , \,\,\,\,\,\,\,\, \lambda_1(t\Omega) = t^{-2+N(1-\frac{2}{q})} \lambda_1(\Omega).
    $$
\end{Lemma}

We may also show the following Faber--Krahn inequality for the first $(2,q)$-eigenvalue, using the Schwarz symmetrisation (see \cite[Section 2.1 and Theorem 3.2.1]{Henrot}).

\begin{prop}\label{FaberKrahn}
    Let $q \in (1,2^*)$, and let $\Omega$ be a bounded open set of $\R^N$. We denote by $\Omega^*$ the ball in $\R^N$ centred at the origin, with measure $|\Omega|$. The first $(2,q)$-eigenvalue satisfies the following Faber--Krahn-type inequality
    $$
        \lambda_1(\Omega) \geq \lambda_1(\Omega^*) = |\Omega|^{1-\frac{2}{q}-\frac{2}{N}} \lambda_1(B_1),
    $$
    where $B_1$ is the ball in $\R^N$ of measure $1$.

    In other words, among all bounded open sets of prescribed measure, the ball minimises the first $(2,q)$-eigenvalue.
\end{prop}

\subsection{The regularity of the $(2,q)$-eigenfunctions}

We finish this section by stating a theorem about the regularity of the solutions of equation \eqref{Eigen2q}.

\begin{prop}\label{Regularity_pq}
    Let $\Omega$ be a bounded open set, and let $q\in (1,2^*)$. Let $u\in H_0^{1}(\Omega)$ be a $(2,q)$-eigenfunction associated with a $(2,q)$-eigenvalue $\lambda$. Therefore, $u\in L^{\infty}(\Omega) \cap C^{1,\alpha}_{loc}(\Omega)$ for some $\alpha \in (0,1)$.
\end{prop}
For the proof see  \cite[Theorem 2.2]{Franzina_Lamberti}.

Since the operator appearing on the left-hand side of equation \eqref{Eigen2q} is simply the Laplacian, we deal with a semilinear equation. In that case, we obtain more regularity, due to the following three estimates.
\begin{prop}[Elliptic estimates]\label{RegularityTheorems}\label{BoundaryRegularityTheorems}
{Let $\Omega\subset\R^N$ be a bounded open set, $1<r<\infty$, $\alpha\in(0,1)$, and let $U\Subset V\Subset\Omega$ be open sets; the constants below depend on $N$, $r$, $\alpha$ and on the sets involved.

(i) \emph{(Calder\'on-Zygmund estimates)} If $u\in H^1_{loc}(\Omega)\cap L^r_{loc}(\Omega)$ and $-\Delta u=f\in L^r_{loc}(\Omega)$, then $u\in W^{2,r}_{loc}(\Omega)$ and $\norm{u}_{W^{2,r}(U)}\leq C\big(\norm{u}_{L^r(V)}+\norm{f}_{L^r(V)}\big)$. Moreover, if $\partial\Omega$ is of class $C^{1,1}$, $u\in H^1_0(\Omega)$ and $f\in L^r(\Omega)$, then $u\in W^{2,r}(\Omega)$ and $\norm{u}_{W^{2,r}(\Omega)}\leq C\norm{f}_{L^r(\Omega)}$.

(ii) \emph{(Morrey's embedding)} If $r>N$, then $W^{2,r}(V)\hookrightarrow C^{1,1-\frac Nr}(\overline U)$, and $W^{2,r}(\Omega)\hookrightarrow C^{1,1-\frac Nr}(\overline\Omega)$ when $\partial\Omega$ is of class $C^{1,1}$.

(iii) \emph{(Schauder estimates)} If $u\in L^\infty_{loc}(\Omega)$ and $-\Delta u=f\in C^{0,\alpha}(\Omega)$, then $u\in C^{2,\alpha}_{loc}(\Omega)$ and $\norm{u}_{C^{2,\alpha}(\overline U)}\leq C\big(\norm{u}_{L^\infty(V)}+\norm{f}_{C^{0,\alpha}(\overline V)}\big)$. Moreover, if $\partial\Omega$ is of class $C^{2,\alpha}$ and $u\in C^0(\overline\Omega)\cap C^2(\Omega)$ solves $-\Delta u=f\in C^{0,\alpha}(\overline\Omega)$ with $u=0$ on $\partial\Omega$, then $u\in C^{2,\alpha}(\overline\Omega)$ and $\norm{u}_{C^{2,\alpha}(\overline\Omega)}\leq C\big(\norm{u}_{C^0(\overline\Omega)}+\norm{f}_{C^{0,\alpha}(\overline\Omega)}\big)$.}
\end{prop}
{For the proofs, see \cite[Theorems 9.11 and 9.19]{Gilbarg_Trudinger} and \cite[Theorems 9.13 and 9.15]{Gilbarg_Trudinger} for (i), \cite[Section 5.6, Theorem 6]{Evans} for (ii), and \cite[Theorem 6.2]{Gilbarg_Trudinger} and \cite[Theorems 6.6 and 6.19]{Gilbarg_Trudinger} for (iii).}

\begin{prop}\label{InteriorRegularity}
    Let $\Omega$ be a bounded open set, and let $q\in (1,2^*)$. Let $u\in H_0^1(\Omega)$ be a $(2,q)$-eigenfunction associated with a $(2,q)$-eigenvalue $\lambda$.
    \begin{itemize}
        \item If $1<q<2$, then $u\in C^{2,q-1}_{loc}(\Omega)$.
        \item If $2\leq q<2^*$, then $u\in C^{2,\alpha}_{loc}(\Omega)$ for all $\alpha<1$.
    \end{itemize}
\end{prop}

\begin{proof}
Let $\lambda\in\sigma$ be a $(2,q)$-eigenvalue. An associated $(2,q)$-eigenfunction $u$ satisfies the following equation
\begin{equation}\label{Eigen2q_proof}
    \left\{
    \begin{array}{cl}
        -\Delta u = \lambda\norm{u}_{L^q(\Omega)}^{2-q}|u|^{q-2}u & \textnormal{ in } \Omega \\
        u = 0 & \textnormal{ on } \partial\Omega.
    \end{array}
    \right.
\end{equation}
We denote by $f$ the second term of \eqref{Eigen2q_proof}: $f = \lambda \norm{u}_{L^q(\Omega)}^{2-q} |u|^{q-2} u$. By Proposition \ref{Regularity_pq},  $u\in L^\infty(\Omega)$. Therefore since $\Omega$ is bounded, then for all $r\in(1,+\infty)$, $f\in L^r(\Omega)$. The interior Calder\'on--Zygmund estimate from Proposition \ref{RegularityTheorems} gives that $u\in W^{2,r}_{loc}(\Omega)$, for all $r\in (1,+\infty)$. By Morrey's embedding, we obtain that, by choosing $r>N$, $u\in C^{1,\alpha}_{loc}(\Omega)$ with $\alpha=1-\frac{N}{r}$. Since $r>N$ is arbitrary,  $u\in C^{1,\alpha}_{loc}(\Omega)$ for all $\alpha<1$.

Denote by $g$ the mapping $g:s\in\R\mapsto|s|^{q-2}s$. We notice that there exists a constant $K=\lambda\norm{u}_{L^q(\Omega)}^{2-q}$ such that $f=K\,g\circ u$. Since the regularity of $f$ depends on that of $g$, which in turn depends on $q$, we distinguish two cases:
\begin{itemize}
    \item If $1<q<2$, $g$ is $(q-1)$-H\"older continuous. Moreover, since $u\in C^{1,\alpha}_{loc}(\Omega)$ for all $\alpha<1$, we have in particular $u\in C^{0,1}_{loc}(\Omega)$. Therefore, $f = K\,g\circ u \in C^{0,q-1}_{loc}(\Omega)$. The interior Schauder estimate from Proposition \ref{RegularityTheorems} implies that $u \in C^{2,q-1}_{loc}(\Omega)$.
    \item If $2\leq q<2^*$, $g$ is more regular; we have $g\in C^{1}(\R)$. Since $u\in C^{1,\alpha}_{loc}(\Omega)$ for all $\alpha<1$, we obtain that $f = K\,g\circ u \in C^{0,\alpha}_{loc}(\Omega)$ for all $\alpha<1$. The interior Schauder estimate from Proposition \ref{RegularityTheorems} implies that $u \in C^{2,\alpha}_{loc}(\Omega)$ for all $\alpha<1$.
\end{itemize}
\end{proof}

\begin{rem}
    By a bootstrap argument, it is also possible to show that $u\in C^\infty_{loc}(\{u\neq 0\} \cap \Omega)$, because the non-linearity $\!g\!:\!s\!\mapsto\!|s|^{q-2}s$ is of class $C^\infty(\R^*)$.
\end{rem}

The global version follows from the global estimates of Proposition \ref{RegularityTheorems}.
\begin{prop}\label{BoundaryRegularity}
    Let $\Omega$ be a bounded open set, and let $q\in (1,2^*)$. Let $u\in H_0^1(\Omega)$ be a $(2,q)$-eigenfunction associated with a $(2,q)$-eigenvalue $\lambda$.
    \begin{itemize}
        \item If $1<q<2$ and $\Omega$ is at least of class $C^{2,q-1}$, then $u\in C^{2,q-1}(\bar{\Omega})$.
        \item If $2\leq q<2^*$ and $\Omega$ is at least of class $C^{2,1}$ then $u\in C^{2,\alpha}(\bar{\Omega})$ for all $\alpha<1$.
    \end{itemize}
\end{prop}

\begin{proof}
    The proof follows exactly the same argument as Proposition \ref{InteriorRegularity}, by replacing local results by their global counterparts given in Proposition \ref{RegularityTheorems}.
\end{proof}

\section{The first non-trivial Fu\v{c}\'\i k curve}

In the following, we restrict our study to the sublinear $q$-eigenvalue problem of the Laplacian, given by \eqref{Eigen2q}. Thus, we fix $q\in(1,2)$, and we consider a bounded open connected set $\Omega\subset\R^N$. We have already defined the first $(2,q)$-eigenvalue $\lambda_1(\Omega)$, via the variational characterisation \eqref{FirstEV}. We would like to define a second $(2,q)$-eigenvalue via a variational characterisation. In that aim, we adopt a method used in \cite{Molle_Passaseo} {based on the study of the Fu\v{c}\'\i k spectrum}.

\begin{defn}
    The \textbf{Fu\v{c}\'\i k $(2,q)$-spectrum} is the set $\Sigma$ of values $(\alpha,\beta)\in\R^2$ such that there exists a weak solution $u\in H_0^1(\Omega)\setminus\{0\}$ of the following problem
    \begin{equation}\label{Fucik2q}
        \left\{
        \begin{array}{cl}
            -\Delta u = - \alpha\norm{u^-}_{L^q(\Omega)}^{2-q}(u^-)^{q-1} + \beta\norm{u^+}_{L^q(\Omega)}^{2-q}(u^+)^{q-1} & \textnormal{ in } \Omega \\
            u = 0 & \textnormal{ on } \partial\Omega,
        \end{array}
        \right.
    \end{equation}
    where $u^+ = \max\{0,u\}$, and $u^-=\max\{0,-u\}$.  The function $u$ satisfying this equation is called a \textbf{$(2,q)$-eigenfunction} associated to the pair $(\alpha,\beta)$.

    In other words, $(\alpha,\beta)\in\Sigma$ if and only if there exists $u\in H_0^1(\Omega)\setminus\{0\}$ such that for all $\varphi\in H_0^1(\Omega)$,
    \begin{equation}\label{Fucik2q_weak}
        \int_\Omega \nabla u \cdot \nabla \varphi \,dx = -\alpha \norm{u^-}_{L^q(\Omega)}^{2-q} \int_\Omega (u^-)^{q-1}\varphi \,dx + \beta \norm{u^+}_{L^q(\Omega)}^{2-q} \int_\Omega (u^+)^{q-1}\varphi \,dx.
    \end{equation}
\end{defn}

If we suppose that $u^+\equiv0$ on $\Omega$, it is trivial to see that $\{\lambda_1(\Omega)\}\times\R$ is contained in the set $\Sigma$, given Propositions \ref{FirstEigenFunctionsArePositive} and \ref{PositiveEigenfunctionsAreFirst}. Similarly, by taking $u^-\equiv0$, we obtain that $\R\times\{\lambda_1(\Omega)\}$ is also contained in $\Sigma$.

Moreover, the set $\Sigma$ is symmetric with respect to the line $\{(\alpha,\beta)\in\R^2: \alpha=\beta\}$. Indeed, if $u$ is an eigenfunction associated to the pair $(\alpha,\beta)$, then $-u$ is associated to the pair $(\beta,\alpha)$.

In addition, if $\alpha\neq \lambda_1(\Omega)$, $\beta\neq \lambda_1(\Omega)$, and $(\alpha,\beta)\in\Sigma$, then $\alpha>\lambda_1(\Omega)$, $\beta>\lambda_1(\Omega)$, and the associated eigenfunction must change sign by Propositions \ref{FirstEigenFunctionsArePositive} and \ref{PositiveEigenfunctionsAreFirst}.

\begin{rem}\label{RescalingEigenfunctions}
    Let $(\alpha,\beta)\in\Sigma$, and let $u\in H_0^1(\Omega)$ be a corresponding eigenfunction. Then, for all $c>0$, the function $cu$ is still an eigenfunction associated to the same pair $(\alpha,\beta)$. Therefore, it is possible to assume that $\norm{u}_{L^q(\Omega)} = 1$, as for a classical eigenvalue problem.

    However, we may not automatically assume that $\norm{u^+}_{L^q(\Omega)} = \norm{u^-}_{L^q(\Omega)}$, since an independent rescaling $k_+u^+-k_-u^-$ of the two nodal parts, with $k_+\neq k_-$, generally creates a jump of the normal derivative across the regular part of the nodal set $\{u=0\}$ of the form $(k_+ - k_-) \partial_\nu u$, and hence an additional singular measure in the distributional Laplacian. This will be important
 when trying to define a second $(2,q)$-eigenfunction via the Fu\v{c}\'\i k $(2,q)$-spectrum.
\end{rem}

The Fu\v{c}\'\i k $(2,q)$-spectrum appears as a generalisation of the classical $(2,q)$-spectrum $\sigma$. In more precise terms, the two following Lemmata hold. Both lemmas follow by writing $|u|^{q-2}u=(u^+)^{q-1}-(u^-)^{q-1}$ and $\norm{u}_{L^q(\Omega)}^q=\norm{u^+}_{L^q(\Omega)}^q+\norm{u^-}_{L^q(\Omega)}^q$ in the two equations; we omit the computations.

\begin{Lemma}\label{EigenvalueToFucik}
    Let $\lambda>\lambda_1(\Omega)$ be a $(2,q)$-eigenvalue, and $u$ be an associated eigenfunction. Therefore, we have
    $$
        \left(\lambda \left(1+\frac{\norm{u^+}_{L^q(\Omega)}^q}{\norm{u^-}_{L^q(\Omega)}^q}\right)^{\frac{2}{q}-1},\lambda \left(1+\frac{\norm{u^-}_{L^q(\Omega)}^q}{\norm{u^+}_{L^q(\Omega)}^q}\right)^{\frac{2}{q}-1}\right) \in \Sigma.
    $$
\end{Lemma}

\begin{Lemma}\label{FucikToEigenvalue}
    Let $(\alpha,\beta)\in \Sigma$, and let $u \in H_0^1(\Omega)$ be an associated sign-changing eigenfunction. If $\alpha\norm{u^-}_{L^q(\Omega)}^{2-q} = \beta\norm{u^+}_{L^q(\Omega)}^{2-q}$, then we can define a classical $(2,q)$-eigenvalue $\lambda \in \sigma$ associated to $u$, given by
    $$
        \lambda = \alpha \left[ 1 + \left( \frac{\alpha}{\beta} \right)^{\frac{q}{2-q}} \right]^{1 - \frac{2}{q}} = \beta \left[ 1 + \left( \frac{\beta}{\alpha} \right)^{\frac{q}{2-q}} \right]^{1 - \frac{2}{q}}.
    $$
\end{Lemma}

\begin{rem}\label{Remark q=2}
     In the case $q=2$:
    \begin{enumerate}
        \item If $\lambda\in \sigma$, then $(\lambda,\lambda) \in \Sigma$.
        \item If $(\alpha,\beta)\in \Sigma$ with $\alpha=\beta=\lambda$, then $\lambda\in\sigma$.
    \end{enumerate}
\end{rem}
As in \cite{Molle_Passaseo}, we will define a second $(2,q)$-eigenvalue through the Fu\v{c}\'\i k $(2,q)$-spectrum. More precisely we will construct a variational characterisation of the smallest (in a sense) pairs $(\alpha,\beta)\in \Sigma$ such that $\alpha\neq\lambda_1(\Omega)$ and $\beta\neq\lambda_1(\Omega)$, called the \textbf{first non-trivial curve} of $\Sigma$. This will be done in defining for all $\beta>\lambda_1(\Omega)$, a family of regularised functionals $(f_{\beta,\varepsilon})_{\varepsilon>0}$ that admit a critical point $u_{\beta,\varepsilon}$ on the set $S=\{u\in H_0^1(\Omega) \, ; \, \int_\Omega (u^-)^q\,dx = 1\}$. By studying the limiting behaviour of the minimiser $u_{\beta,\varepsilon}$ when $\varepsilon$ goes to $0$, we will obtain a variational characterisation for the first curve of $\Sigma$. However, to define a second $(2,q)$-eigenvalue via this procedure, we will need an additional nodal-balance condition.

\subsection{Defining and minimising the family $(f_{\beta,\varepsilon})_\varepsilon$}\label{Section First Curve: Defining The Family}

For all $\beta>0$ and all $\varepsilon>0$, we define the functional $f_{\beta,\varepsilon}: H_0^1(\Omega)\rightarrow\R$ as
\begin{equation}
    f_{\beta,\varepsilon}(u) = \int_\Omega |\nabla u|^2\,dx - q \left[\int_\Omega G_{\beta,\varepsilon}(u)\,dx\right]^{\frac{2}{q}},
\end{equation}
where  $G_{\beta,\varepsilon}(t) = \int_0^t g_{\beta,\varepsilon}(\tau)\,d\tau$, with
$$
    g_{\beta,\varepsilon}(\tau) = \left\{
    \begin{array}{cl}
        0, & \textnormal{ if } \tau\leq\varepsilon \\
        \beta^{\frac{q}{2}}q^{1-\frac{q}{2}} \left(\tau^{q-1} - \varepsilon^{q-1}\right), & \textnormal{ if } \tau>\varepsilon.
    \end{array}
    \right.
$$
By simple computation, we have for all $u\in H_0^1(\Omega)$ and all $\varphi\in H_0^1(\Omega)$,
$$
    f_{\beta,\varepsilon}'(u)[\varphi] = 2\int_\Omega \nabla u\cdot\nabla\varphi\,dx - 2\left[\int_\Omega G_{\beta,\varepsilon}(u)\,dx\right]^{\frac{2}{q}-1}\int_\Omega g_{\beta,\varepsilon}(u)\varphi\,dx.
$$
Finally, we consider the set
$$
S = \{u\in H_0^1(\Omega) ; \int_\Omega (u^-)^q\,dx=1\}
$$
and the subset $M^{\beta,\varepsilon}$ of functions $u$ in $S$ such that $u^+\not\equiv 0$ and $f_{\beta,\varepsilon}'(u)[u^+]=0$.

\begin{rem}\label{Properties of G and g}
    By the definition of $g_{\beta,\varepsilon}$, we note that for all $\varepsilon>0$, all $\beta>0$, and all functions $v\in H_0^1(\Omega)$,
    $$
        g_{\beta,\varepsilon}(v) = g_{\beta,\varepsilon}(v^+-v^-) = g_{\beta,\varepsilon}(v^+) \,\,;
    $$
    $$
         G_{\beta,\varepsilon}(v) = \int_0^{v^+-v^-} g_{\beta,\varepsilon}(\tau)\,d\tau = \int_0^{v^+} g_{\beta,\varepsilon}(\tau)\,d\tau = G_{\beta,\varepsilon}(v^+)
    $$
    (since the supports of $v^+$ and $v^-$ are disjoint). Similarly, for all $t>0$, we will have
    $$
        g_{\beta,\varepsilon}(tv^+ - v^-) = g_{\beta,\varepsilon}(tv^+) \,\,\,\, ; \,\,\,\,       G_{\beta,\varepsilon}(tv^+ - v^-) = G_{\beta,\varepsilon}(tv^+).
    $$
\end{rem}

\begin{Lemma}\label{Lemma1}
    Let $\varepsilon>0$ and $\beta>0$. If $u\in S$ is a sign-changing function which is a critical point of $f_{\beta,\varepsilon}$ on $S$, then $u\in M^{\beta,\varepsilon}$. Moreover, for all $u\in M^{\beta,\varepsilon}$, the function $t\mapsto f_{\beta,\varepsilon}'(-u^- + tu^+)[u^+]$ is positive for $t\in(0,1)$, and negative for $t>1$.
\end{Lemma}
\begin{proof}
For $s\geq0$ we have $-u^-+su^+\in S$, so that $s\mapsto f_{\beta,\varepsilon}(-u^-+su^+)$ is differentiable with a critical point at $s=1$, that is $f'_{\beta,\varepsilon}(u)[u^+]=0$; since $u^+\not\equiv0$, $u\in M^{\beta,\varepsilon}$.

Let now $u\in M^{\beta,\varepsilon}$. Then $|\{u^+>\varepsilon\}|>0$: otherwise $g_{\beta,\varepsilon}(u^+)\equiv0$ and $f'_{\beta,\varepsilon}(u)[u^+]=0$ would give $\int_\Omega|\nabla u^+|^2\,dx=0$. The map $\tau\mapsto\tau^{1-q}g_{\beta,\varepsilon}(\tau)=\beta^{\frac q2}q^{1-\frac q2}\big(1-(\frac\varepsilon\tau)^{q-1}\big)$ is non-decreasing on $(0,+\infty)$ and strictly increasing on $(\varepsilon,+\infty)$. Hence, for $0<t<1$, pointwise in $\Omega$,
$$
g_{\beta,\varepsilon}(tu^+)\leq t^{q-1}g_{\beta,\varepsilon}(u^+),\qquad G_{\beta,\varepsilon}(tu^+)=t\int_0^{u^+}g_{\beta,\varepsilon}(t\sigma)\,d\sigma\leq t^{q}G_{\beta,\varepsilon}(u^+),
$$
with strict inequalities on $\{tu^+>\varepsilon\}$. Since $\int_\Omega G_{\beta,\varepsilon}(u^+)\,dx>0$ and $\int_\Omega g_{\beta,\varepsilon}(u^+)u^+\,dx>0$, integrating and multiplying we obtain
$$
\Big[\int_\Omega G_{\beta,\varepsilon}(tu^+)\,dx\Big]^{\frac2q-1}\int_\Omega g_{\beta,\varepsilon}(tu^+)u^+\,dx<t\Big[\int_\Omega G_{\beta,\varepsilon}(u^+)\,dx\Big]^{\frac2q-1}\int_\Omega g_{\beta,\varepsilon}(u^+)u^+\,dx
$$
(the left-hand side vanishes if $|\{tu^+>\varepsilon\}|=0$). By Remark \ref{Properties of G and g}, this means
\[
\begin{split}
    f'_{\beta,\varepsilon}(-u^-+tu^+)[u^+]& = 2t\int_\Omega|\nabla u^+|^2\,dx-2\Big[\int_\Omega G_{\beta,\varepsilon}(tu^+)\,dx\Big]^{\frac2q-1}\int_\Omega g_{\beta,\varepsilon}(tu^+)u^+\,dx\\&>t\,f'_{\beta,\varepsilon}(u)[u^+]=0.
\end{split}
\]
For $t>1$ the two pointwise inequalities are reversed, strictly on $\{u^+>\varepsilon\}$, and the same computation gives $f'_{\beta,\varepsilon}(-u^-+tu^+)[u^+]<0$.
\end{proof}

\begin{cor}\label{cor:maxpositive}
Let $\varepsilon,\beta>0$ and $u\in M^{\beta,\varepsilon}$. Then $t\mapsto f_{\beta,\varepsilon}(tu^+)$ attains its maximum on $[0,+\infty)$ exactly at $t=1$, and $f_{\beta,\varepsilon}(u^+)>0$. In particular, $u$ is the only maximum point of $f_{\beta,\varepsilon}$ on the set $\{-u^-+tu^+ ; t>0\}\subset S$.
\end{cor}

\begin{proof}
Since $u^+u^-=0$, $F(t)=f_{\beta,\varepsilon}(-u^-+tu^+)=f_{\beta,\varepsilon}(-u^-)+f_{\beta,\varepsilon}(tu^+)$, and $F'(t)=f'_{\beta,\varepsilon}(-u^-+tu^+)[u^+]$ is positive on $(0,1)$ and negative on $(1,+\infty)$ by Lemma \ref{Lemma1}. Hence $f_{\beta,\varepsilon}(u^+)>f_{\beta,\varepsilon}(0)=0$.
\end{proof}

Before proving the next proposition, we prove a lemma related to the limiting behaviour of the first $(2,q)$-eigenvalue on some sequences of bounded open sets.

\begin{Lemma}\label{DecreasingSequenceOfLambda1}
    Let $x_0\in\partial\Omega$. For all $r\in (0,\frac{1}{2} diam(\Omega))$, let $\Omega_r=\Omega\setminus \bar{B}(x_0,r)$. The family $\lambda_1(\Omega_r)$ converges to $\lambda_1(\Omega)$ as $r$ goes to $0$.
\end{Lemma}

\begin{proof}
    For all $r_1$ and $r_2$ such that $r_2 > r_1 > 0$, since $\Omega_{r_2} \subseteq \Omega_{r_1} \subseteq \Omega$, Corollary \ref{StrictMonotonicityEV1} yields $\lambda_1(\Omega_{r_2}) \geq \lambda_1(\Omega_{r_1}) \geq \lambda_1(\Omega)$.  Thus, the function $r\mapsto\lambda_1(\Omega_r)$ is non-increasing as $r$ decreases to $0$, so the corresponding limit exists.

    Conversely, let $\delta>0$, and let $\phi\in C_c^\infty(\Omega)$ whose Rayleigh quotient is smaller than $\lambda_1(\Omega)+\delta$. Since $Supp(\phi) \Subset \Omega$ and $x_0\in\partial\Omega$, there exists $r_\phi>0$ such that for all $r\in(0,r_\phi)$, we have $Supp(\phi) \subset \Omega_r$. Thus, for all $r\in(0,r_\phi)$, $\phi\in H_0^1(\Omega_r)$, and formula \eqref{FirstEV} implies that
    $\displaystyle \lambda_1(\Omega_r) \leq \frac{\int_\Omega |\nabla \phi|^2 \,dx}{\norm{\phi}_{L^q(\Omega)}^{2}} \leq \lambda_1(\Omega)+\delta$.
    Hence, $\limsup_{r\rightarrow 0}\lambda_1(\Omega_r) \leq \lambda_1(\Omega) + \delta$.

    Thus, we conclude that $\lambda_1(\Omega_r)$ converges to $\lambda_1(\Omega)$ as $r$ goes to $0$.
\end{proof}

\begin{prop}\label{M_nonempty}
    For all $\beta>0$ and all $\varepsilon>0$,
    $
        M^{\beta,\varepsilon}\neq\emptyset \,\, \Longleftrightarrow \,\, \beta>\lambda_1(\Omega).
    $
\end{prop}

\begin{proof}
We first prove that $M^{\beta,\varepsilon}\neq\emptyset$ implies $\beta>\lambda_1(\Omega)$. This will be done by showing that $M^{\beta,\varepsilon}=\emptyset$ if $\beta \leq \lambda_1(\Omega)$.
For all $\beta>0$, $\varepsilon>0$, and $\tau\geq 0$,  $g_{\beta,\varepsilon}(\tau)\leq \beta^{\frac{q}{2}}q^{1-\frac{q}{2}} \tau^{q-1}$. Let $u\in S$ be a sign-changing function (i.e. such that $u^+\not\equiv 0$). Thus,
$$
    G_{\beta,\varepsilon}(u^+) = \int_0^{u^+} g_{\beta,\varepsilon}(\tau)\,d\tau \leq \beta^{\frac{q}{2}}q^{1-\frac{q}{2}} \int_0^{u^+} \tau^{q-1}\,d\tau = \beta^{\frac{q}{2}}q^{-\frac{q}{2}} (u^+)^q.
$$
We then obtain
\[
\begin{split}
    \left[\int_\Omega G_{\beta,\varepsilon}(u^+)\,dx\right]^{\frac{2}{q}-1}\int_\Omega g_{\beta,\varepsilon}(u^+)u^+\,dx
    &\leq \left[\beta^{\frac{q}{2}}q^{-\frac{q}{2}} \int_\Omega (u^+)^q\,dx\right]^{\frac{2}{q}-1} \beta^{\frac{q}{2}}q^{1-\frac{q}{2}}\int_\Omega (u^+)^{q-1}u^+\,dx \\
    &=\beta \left[\int_\Omega (u^+)^q\,dx\right]^{\frac{2}{q}}.
\end{split}
\]
By the variational characterisation \eqref{FirstEV},  $\int_\Omega |\nabla u^+|^2\,dx \geq \lambda_1(\Omega) \left[\int_\Omega (u^+)^q\,dx\right]^{\frac{2}{q}}$. The inequality is strict. Indeed, equality would make $u^+$ a first eigenfunction on $\Omega$; by Proposition \ref{FirstEigenFunctionsArePositive}, it would be strictly positive in the connected domain $\Omega$, whereas $u^+ \equiv 0$ in the set $\{u^->0\}$, which has positive measure because $\norm{u^-}_{L^q(\Omega)} = 1$.

Therefore, if $\beta\leq\lambda_1(\Omega)$, we have
\[
\begin{split}    \int_\Omega |\nabla u^+|^2\,dx &> \lambda_1(\Omega) \left[\int_\Omega (u^+)^q\,dx\right]^{\frac{2}{q}} \geq \beta \left[\int_\Omega (u^+)^q\,dx\right]^{\frac{2}{q}} \\ &\geq \left[\int_\Omega G_{\beta,\varepsilon}(u^+)\,dx\right]^{\frac{2}{q}-1}\int_\Omega g_{\beta,\varepsilon}(u^+)u^+\,dx.\end{split}
\]
Thus, $f_{\beta,\varepsilon}'(u)[u^+]>0$, which shows that $M^{\beta,\varepsilon} = \emptyset$ for all $\varepsilon>0$ and all $\beta\leq \lambda_1(\Omega)$.

We now prove that $M^{\beta,\varepsilon}\neq\emptyset$ if $\beta>\lambda_1(\Omega)$.
 To do that we construct a function $\bar{u}\in S$ such that
$$
    \lim_{t\rightarrow 0} \frac{1}{t^2}f_{\beta,\varepsilon}(t\bar{u}^+) >0\,,\quad
    \lim_{t\rightarrow +\infty} \frac{1}{t^2}f_{\beta,\varepsilon}(t\bar{u}^+) <0\,.
$$
Therefore, for all $\varepsilon>0$,
there exists $\bar{t}>0$ such that for all $t\geq 0$, $f_{\beta,\varepsilon}(-\bar{u}^-+\bar{t}\bar{u}^+) \geq f_{\beta,\varepsilon}(-\bar{u}^-+t\bar{u}^+)$. By Lemma \ref{Lemma1},  $-\bar{u}^-+\bar{t}\bar{u}^+ \in M^{\beta,\varepsilon}$.
\\
Under the assumption
 $\beta>\lambda_1(\Omega)$,  we are able to construct a function $\bar{u}\in S$ such that $\int_\Omega |\nabla \bar{u}^+|^2\,dx < \beta \left[\int_\Omega (\bar{u}^+)^q\,dx\right]^{\frac{2}{q}}$, as follows. Let $x_0\in\partial\Omega$. For all $r\in (0,\frac{1}{2} diam(\Omega))$, we set $\Omega_r=\Omega\setminus \bar{B}(x_0,r)$ and $\tilde{\Omega}_r=\Omega\cap B(x_0,r)$. We consider the eigenfunctions $e_1(\Omega_r)$ and $e_1(\tilde{\Omega}_r)$, normalised in $L^q$, associated with $\lambda_1(\Omega_r)$ and $\lambda_1(\tilde{\Omega}_r)$, respectively. Replacing these first eigenfunctions by their absolute values, we may assume that they are non-negative on their respective domains, $\Omega_r$ and $\tilde{\Omega}_r$.
 (We may notice that no connectedness of $\Omega_r$ or $\tilde{\Omega}_r$ is needed here.) We then define the following function
$$
    \bar{u} = \left\{
    \begin{array}{cl}
        e_1(\Omega_r) & \textnormal{ on } \Omega_r \\
        -e_1(\tilde{\Omega}_r) & \textnormal{ on } \tilde{\Omega}_r.
    \end{array}
    \right.
$$
Therefore, we have $\bar{u}\in S$, and
$$
    \int_\Omega |\nabla \bar{u}^+|^2\,dx = \int_{\Omega_r} |\nabla e_1(\Omega_r)|^2\,dx = \lambda_1(\Omega_r) \left[\int_{\Omega_r} (e_1(\Omega_r))^q\,dx\right]^{\frac{2}{q}} = \lambda_1(\Omega_r).
$$
By Lemma \ref{DecreasingSequenceOfLambda1}, $\lambda_1(\Omega_r)$ decreases and tends to $\lambda_1(\Omega)$ when $r$ goes towards $0$; thus, we may choose $r$ small enough such that $\beta>\lambda_1(\Omega_r)$, since $\beta>\lambda_1(\Omega)$. This choice of $r$ gives a suitable function $\bar{u}\in S$.
\\
For all $t>0$,
\[
\begin{split}
    \frac{1}{t^2}f_{\beta,\varepsilon}(t\bar{u}^+)
    &= \int_\Omega |\nabla \bar{u}^+|^2\,dx - \frac{\beta q^{\frac{2}{q}}}{t^2} \left[\int_\Omega \left(\frac{1}{q}t^q (\bar{u}^+)^q - \varepsilon^{q-1}t\bar{u}^+ + \left(1-\frac{1}{q}\right) \varepsilon^{q}\right)  \,\,\,\, \chi_{\{t\bar{u}^+>\varepsilon\}} \,dx\right]^{\frac{2}{q}} \\
    &= \int_\Omega |\nabla \bar{u}^+|^2\,dx - \beta \left[\int_\Omega \left((\bar{u}^+)^q - \frac{q\varepsilon^{q-1}}{t^{q-1}}\bar{u}^+ + \frac{q-1}{t^q} \varepsilon^{q}\right)  \,\,\,\, \chi_{\{t\bar{u}^+>\varepsilon\}} \,dx\right]^{\frac{2}{q}}.
\end{split}
\]
Since $q>1$, we have
$$
    \lim_{t\rightarrow +\infty} \frac{1}{t^2}f_{\beta,\varepsilon}(t\bar{u}^+) = \int_\Omega |\nabla \bar{u}^+|^2\,dx - \beta \left[\int_\Omega (\bar{u}^+)^q\,dx\right]^{\frac{2}{q}} < 0.
$$

Since the supports of $\bar{u}^+$ and $\bar{u}^-$ are disjoint, then for all $t\geq 0$, $f_{\beta,\varepsilon}(-\bar{u}^-+t\bar{u}^+) = f_{\beta,\varepsilon}(-\bar{u}^-) + f_{\beta,\varepsilon}(t\bar{u}^+)$. We then have
\[
\begin{split}
    f_{\beta,\varepsilon}(t\bar{u}^+)
    &= t^2\int_\Omega |\nabla \bar{u}^+|^2\,dx - q \left[\int_\Omega G_{\beta,\varepsilon}(t\bar{u}^+)\,dx\right]^{\frac{2}{q}} \\
    &= t^2\int_\Omega |\nabla \bar{u}^+|^2\,dx - q \left[\int_\Omega \int_\varepsilon^{t\bar{u}^+} \beta^{\frac{q}{2}}q^{1-\frac{q}{2}}\left(\tau^{q-1} - \varepsilon^{q-1}\right) \,d\tau \,\,\,\, \chi_{\{t\bar{u}^+>\varepsilon\}} \,dx\right]^{\frac{2}{q}}.
\end{split}
\]
Since $\varepsilon>0$ and $e_1(\Omega_r)\in L^\infty (\Omega)$ (by Proposition \ref{Regularity_pq}),  we have $G_{\beta,\varepsilon}(t \bar{u}^+) \equiv 0$ for $0<t<\frac{\varepsilon}{\norm{\bar{u}^+}_{L^\infty(\Omega)}}$. Hence,
$$
    \lim_{t\rightarrow 0} \frac{1}{t^2}f_{\beta,\varepsilon}(t\bar{u}^+) = \int_\Omega |\nabla \bar{u}^+|^2\,dx >0.
$$
\end{proof}

\begin{rem}\label{rem:testfunction}
    For every $\beta>\lambda_1(\Omega)$, the function $\bar u$ constructed in the second part of the previous proof satisfies $\norm{\bar u^+}_{L^q(\Omega)}=\norm{\bar u^-}_{L^q(\Omega)}=1$ and $\int_\Omega|\nabla\bar u^+|^2\,dx=\lambda_1(\Omega_r)<\beta$. It will be used again in the proof of Proposition \ref{ConvergenceSequence}.
\end{rem}

Having now shown that the set $M^{\beta,\varepsilon}$ is non-empty for all $\beta>\lambda_1(\Omega)$ and all $\varepsilon>0$, we are interested in minimising $f_{\beta,\varepsilon}$ on that set.

\begin{Lemma}\label{IntermediateValueTheorem_f}
    Let $\beta>0$ and $\varepsilon>0$. Let $w\in H_0^1(\Omega)$ such that $w\geq 0$, $w \not\equiv 0$, and $\int_\Omega |\nabla w|^2 \,dx < \beta \norm{w}_{L^q(\Omega)}^2$. Then, there exists a unique $t_\varepsilon>0$ such that $f'_{\beta,\varepsilon}(t_\varepsilon w)[w]=0$. Moreover, $t\mapsto f_{\beta,\varepsilon}(tw)$ attains its maximum exactly at $t_\varepsilon$, and $\lim_{\varepsilon\to0}f_{\beta,\varepsilon}(t_\varepsilon w) = 0$.
\end{Lemma}

\begin{proof}
    We define the mapping $\Phi: t\in(0,+\infty) \, \mapsto f_{\beta,\varepsilon} (tw)$. Its derivative is, for all $t>0$,
    $$
        \Phi'(t) = f_{\beta,\varepsilon}'(tw)[w] = 2t \int_\Omega |\nabla w|^2 \,dx - 2\left[\int_\Omega G_{\beta,\varepsilon}(tw)\,dx\right]^{\frac{2}{q}-1}\int_\Omega g_{\beta,\varepsilon}(tw)w\,dx.
    $$
For $t>0$, the bounds $0\leq g_{\beta,\varepsilon}(\tau)\leq\beta^{\frac q2}q^{1-\frac q2}\tau^{q-1}$, $0\leq G_{\beta,\varepsilon}(\tau)\leq\beta^{\frac q2}q^{-\frac q2}\tau^{q}$ and $g_{\beta,\varepsilon}(tw)=G_{\beta,\varepsilon}(tw)=0$ on $\{tw\leq\varepsilon\}$ give
$$
\frac1t\Big[\int_\Omega G_{\beta,\varepsilon}(tw)\,dx\Big]^{\frac2q-1}\int_\Omega g_{\beta,\varepsilon}(tw)w\,dx\leq\beta\Big[\int_{\{tw>\varepsilon\}}w^q\,dx\Big]^{\frac2q}\xrightarrow[t\to0]{}0,
$$
because $|\{tw>\varepsilon\}|\to0$ as $t\to0$ and $w\in L^q(\Omega)$. Hence $\frac{1}{2t}\Phi'(t)\to\int_\Omega|\nabla w|^2\,dx>0$, and $\Phi'(t)>0$ for $t>0$ small enough.

On the other hand, for all $t>0$,
\[
\begin{split}
    \frac{1}{2t}\Phi'(t)&= \int_\Omega |\nabla w|^2\,dx - \frac{\beta}{t} \Big[\int_\Omega \big(t^q w^q - q\varepsilon^{q-1}tw + (q-1) \varepsilon^{q}\big) \chi_{\{tw>\varepsilon\}} \,dx\Big]^{\frac{2}{q}-1}
    \\ &\qquad\times \int_\Omega \big( t^{q-1}w^q - \varepsilon^{q-1}w \big) \chi_{\{tw>\varepsilon\}} \,dx
    \\ &= \int_\Omega |\nabla w|^2\,dx - \beta \Big[\int_\Omega \Big(w^q - \frac{q\varepsilon^{q-1}}{t^{q-1}}w + \frac{(q-1)}{t^{q}} \varepsilon^{q}\Big) \chi_{\{tw>\varepsilon\}} \,dx\Big]^{\frac{2}{q}-1}
    \\ &\qquad\times \int_\Omega \Big( w^q - \frac{\varepsilon^{q-1}}{t^{q-1}}w \Big) \chi_{\{tw>\varepsilon\}} \,dx\,.
\end{split}
\]
Therefore, since $q>1$, we obtain
    $$
        \lim_{t\rightarrow +\infty} \frac{1}{2t}\Phi'(t) = \int_\Omega |\nabla w|^2\,dx - \beta \left[\int_\Omega w^q \,dx\right]^{\frac{2}{q}} < 0.
    $$
    Hence, for $t>0$ small enough, $\Phi'(t)>0$, and for $t$ large enough, $\Phi'(t)<0$. Since $\Phi'$ is continuous, the intermediate value theorem gives the existence of a real number $t_\varepsilon>0$ such that $\Phi'(t_\varepsilon)=0$.

We now prove the uniqueness. We have
$$
    \frac1{2t}\Phi'(t)=\int_\Omega|\nabla w|^2\,dx-H(t),
$$
where
$$
    H(t)=\frac1t\Big[\int_\Omega G_{\beta,\varepsilon}(tw)\,dx\Big]^{\frac2q-1}\int_\Omega g_{\beta,\varepsilon}(tw)w\,dx
$$
is non-decreasing on $(0,+\infty)$ and strictly increasing on $\{t:|\{tw>\varepsilon\}|>0\}$, by the monotonicity of $\tau\mapsto\tau^{1-q}g_{\beta,\varepsilon}(\tau)$ as in the proof of Lemma \ref{Lemma1}; since $\Phi'(t)>0$ where $|\{tw>\varepsilon\}|=0$, the zero $t_\varepsilon$ of $\Phi'$ is unique.

Hence, for $0<t<t_\varepsilon$,  $f_{\beta,\varepsilon}'(tw)[w]>0$, and for $t>t_\varepsilon$, $f_{\beta,\varepsilon}'(tw)[w]<0$. Therefore, $f_{\beta,\varepsilon}$ attains its maximum exactly at $t_\varepsilon$.

It remains to show that $\lim_{\varepsilon\to0}f_{\beta,\varepsilon}(t_\varepsilon w) = 0$. For every $s\geq 0$, we have $G_{\beta,\varepsilon}(s)\ge \beta^{\frac q2}q^{-\frac q2}(s^q-q\varepsilon^{q-1}s)$. Hence, there exists a constant $C>0$ such that for every $t\geq 0$,
$$
    f_{\beta,\varepsilon}(t w)\leq t^2 \Bigl(\int_\Omega|\nabla w|^2\,dx-\beta\norm{w}_{L^q(\Omega)}^2 \Bigr) + C\beta\,\varepsilon^{q-1} t^{3-q}\norm{w}_{L^1(\Omega)},
$$
whenever the lower bound for $\int_\Omega G_{\beta,\varepsilon}(t w)\,dx$ is non-negative, while $f_{\beta,\varepsilon}(t w)\leq t^2\int_\Omega|\nabla w|^2\,dx$ always. Since $\int_\Omega |\nabla w|^2\,dx < \beta \norm{w}_{L^q(\Omega)}^2$ and $3-q<2$, maximising over $t$ yields $\max_{t\geq 0}f_{\beta,\varepsilon}(t w)\leq C\varepsilon^{2}$, and the claim follows.
\end{proof}

\begin{Lemma}\label{lem:boundedminimising}
    Let $\beta>\lambda_1(\Omega)$ and $\varepsilon>0$, and let $(u_n)_{n\in\N}$ be a minimising sequence for $f_{\beta,\varepsilon}$ on $M^{\beta,\varepsilon}$. Then $(u_n^-)_n$ and $(u_n^+)_n$ are bounded in $H_0^1(\Omega)$.
\end{Lemma}

\begin{proof}
Let $(u_n)_{n\in\N}$ be a minimising sequence in $M^{\beta,\varepsilon}$. For all $n\in\N$, $f_{\beta,\varepsilon}(u_n) = f_{\beta,\varepsilon}(-u_n^-)+f_{\beta,\varepsilon}(u_n^+)$, with $f_{\beta,\varepsilon}(-u_n^-) = \int_\Omega |\nabla u_n^-|^2\,dx \geq \lambda_1(\Omega) \left[\int_\Omega (u_n^-)^q\,dx\right]^\frac{2}{q} = \lambda_1(\Omega)$.
By Corollary \ref{cor:maxpositive}, $f_{\beta,\varepsilon}(u_n^+)=\max\{f_{\beta,\varepsilon}(tu_n^+)\,;\,t\geq0\}>0$.

Therefore, by considering that $\sup\{f_{\beta,\varepsilon}(u_n) ; n\in\N\}<+\infty$, we obtain that
\begin{equation}\label{In-}
    \lambda_1(\Omega)\leq \liminf_{n\rightarrow+\infty} f_{\beta,\varepsilon}(-u_n^-) \leq \limsup_{n\rightarrow+\infty} f_{\beta,\varepsilon}(-u_n^-) < +\infty,
\end{equation}
and
\begin{equation}\label{In+}
    0 \leq \liminf_{n\rightarrow+\infty} f_{\beta,\varepsilon}(u_n^+) \leq \limsup_{n\rightarrow+\infty} f_{\beta,\varepsilon}(u_n^+) < +\infty.
\end{equation}
Since for all $n\in\N$, $u_n\in H_0^1(\Omega)$ and $f_{\beta,\varepsilon}(-u_n^-) = \int_\Omega |\nabla u_n^-|^2\,dx$, inequality \eqref{In-} and Poincar\'e's inequality imply that the sequence $(u_n^-)_n$ is bounded $H_0^1(\Omega)$.

We must now show that $(u_n^+)_n$ is also bounded in $H_0^1(\Omega)$. Since for all $n\in\N$, $u_n\in M^{\beta
,\varepsilon}$, we have $f_{\beta,\varepsilon}'(u_n)[u_n^+]=0$. Therefore
\begin{equation}\label{In1}
    \int_\Omega |\nabla u_n^+|^2\,dx = \left[\int_\Omega G_{\beta,\varepsilon}(u_n^+)\,dx\right]^{\frac{2}{q}-1}\int_\Omega g_{\beta,\varepsilon}(u_n^+)u_n^+\,dx \leq \beta \left[\int_\Omega (u_n^+)^q \,dx\right]^{\frac{2}{q}}.
\end{equation}
It suffices to show that $(u_n^+)_n$ is bounded in $L^q(\Omega)$. By contradiction, we assume that, up to extraction of a subsequence, $\lim_{n\rightarrow+\infty} \int_\Omega (u_n^+)^q\,dx=+\infty$. We set for all $n\in\N$, $v_n = \frac{u_n^+}{\norm{u_n^+}_{L^q(\Omega)}}$. Inequality \eqref{In1} implies that for all $n\in\N$, $\int_\Omega |\nabla v_n|^2\,dx \leq \beta$. Therefore, up to another extraction, $(v_n)_n$ converges to a function $v\in H_0^1(\Omega)$, strongly in $L^q(\Omega)$, weakly in $H_0^1(\Omega)$, and a.e. in $\Omega$. Consequently, $\int_\Omega |\nabla v|^2\,dx \leq \beta$, $\int_\Omega v^q\,dx=1$ and $v\geq 0$ in $\Omega$. Since $f_{\beta,\varepsilon}'(u_n)[u_n^+]=0$ for all $n\in\N$, and $\lim_{n\rightarrow+\infty} \int_\Omega (u_n^+)^q\,dx=+\infty$, we have
\[
\begin{split}
    \int_\Omega |\nabla v_n|^2\,dx
    = \frac{\left[\int_\Omega G_{\beta,\varepsilon}(u_n^+)\,dx\right]^{\frac{2}{q}-1}\int_\Omega g_{\beta,\varepsilon}(u_n^+)u_n^+\,dx}{\left[\int_\Omega (u_n^+)^q\,dx\right]^{\frac{2}{q}}}
    \xrightarrow[n\rightarrow +\infty]{} \beta.
\end{split}
\]
In addition, since $(v_n)_n$ converges to $v$ strongly in $L^q(\Omega)$ and a.e. in $\Omega$, the dominated convergence theorem yields that for all $t\geq 0$,
\[
\begin{split}
    \lim_{n\rightarrow+\infty} f_{\beta,\varepsilon}'(tv_n)[v_n]
    &= \lim_{n\rightarrow+\infty} \left[2t\int_\Omega |\nabla v_n^+|^2\,dx - 2\left[\int_\Omega G_{\beta,\varepsilon}(tv_n)\,dx\right]^{\frac{2}{q}-1}\int_\Omega g_{\beta,\varepsilon}(tv_n)v_n\,dx\right] \\
    &= 2t\beta - 2\left[\int_\Omega G_{\beta,\varepsilon}(tv)\,dx\right]^{\frac{2}{q}-1}\int_\Omega g_{\beta,\varepsilon}(tv)v\,dx \\
    &= 2\int_\Omega \left[t\beta v^{q-1} - \left[\int_\Omega G_{\beta,\varepsilon}(tv)\,dx\right]^{\frac{2}{q}-1} g_{\beta,\varepsilon}(tv)\right]v\,dx,
\end{split}
\]
because $\int_\Omega v^q\,dx=1$. Note that, by the bounds $G_{\beta,\varepsilon}(s)\leq\beta^{\frac q2}q^{-\frac q2}s^q$
and $g_{\beta,\varepsilon}(s)\leq\beta^{\frac q2}q^{1-\frac q2}s^{q-1}$ for all $s\geq 0$, together
with $\int_\Omega v^q\,dx=1$, we have
$$
    \Bigl[\int_\Omega G_{\beta,\varepsilon}(tv)\,dx\Bigr]^{\frac 2q-1} g_{\beta,\varepsilon}(tv)\;\leq\;\beta\,t\,v^{q-1} \qquad\text{a.e. in }\Omega,\ \forall t>0 .
$$
Moreover, since for $s>\varepsilon$, $g_{\beta,\varepsilon}(s)=\beta^{\frac q2}q^{1-\frac q2}s^{q-1}\bigl(1 - \left(\frac{\varepsilon}{s}\right)^{q-1}\bigr)$ and $\int_\Omega G_{\beta,\varepsilon}(tv)\,dx = \beta^{\frac q2}q^{-\frac q2}t^q\bigl(1+O(t^{1-q})\bigr)$ as $t\to+\infty$, a direct computation shows that
$$
    h(t):=2t\beta-2\Bigl[\int_\Omega G_{\beta,\varepsilon}(tv)\,dx\Bigr]^{\frac 2q-1} \int_\Omega g_{\beta,\varepsilon}(tv)\,v\,dx \;\geq\; c\,t^{2-q}\quad\text{for $t$ large},
$$
for a suitable $c=c(\beta,q,\varepsilon,v)>0$; in particular $\int_0^\tau h(t)\,dt\to+\infty$ as $\tau\to+\infty$, because $q<2$. For all $n\in\N$, set $t_n = \norm{u_n^+}_{L^q(\Omega)}$. As in the proof of Lemma \ref{Lemma1}, we have that $f_{\beta,\varepsilon}'(tv_n)[v_n]>0$, for all $t\in(0,t_n[$. Since $t_n\rightarrow+\infty$, and by the fundamental theorem of calculus, we obtain
\[
\begin{split}
    \liminf_{n\rightarrow+\infty} f_{\beta,\varepsilon}(u_n^+)
    &= \liminf_{n\rightarrow+\infty} \int_0^{t_n} f_{\beta,\varepsilon}'(tv_n)[v_n]\,dt \\
    &\geq \int_0^\tau \left[2t\beta - 2\left[\int_\Omega G_{\beta,\varepsilon}(tv)\,dx\right]^{\frac{2}{q}-1} \int_\Omega g_{\beta,\varepsilon}(tv)v\,dx\right]\,dt \,\,,\,\, \forall\tau>0
\end{split}
\]
Therefore, when $\tau\rightarrow+\infty$, we must have $\lim_{n\rightarrow+\infty} f_{\beta,\varepsilon}(u_n^+)=+\infty$. This contradicts inequality \eqref{In+}.
Hence, the sequence $(u_n^+)_n$ is bounded in $H_0^1(\Omega)$.
\end{proof}

\begin{prop}\label{InfAchieved}
    For all $\beta>\lambda_1(\Omega)$, the infimum of $f_{\beta,\varepsilon}$ on $M^{\beta,\varepsilon}$ is achieved for all $\varepsilon>0$.
\end{prop}

\begin{proof}
We fix $\beta>\lambda_1(\Omega)$ and $\varepsilon>0$. Let $(u_n)_{n\in\N}$ be a minimising sequence in $M^{\beta,\varepsilon}$. By Lemma \ref{lem:boundedminimising}, $(u_n)_n$ is bounded in $H_0^1(\Omega)$. Consequently, up to extraction of a subsequence, $(u_n)_n$ converges to a function $u\in H_0^1(\Omega)$, strongly in $L^q(\Omega)$, weakly in $H_0^1(\Omega)$, and a.e. in $\Omega$. Therefore, we have $\int_\Omega (u^-)^q\,dx = 1$. It remains to show that $u^+\not\equiv 0$. By contradiction, we assume that $u^+\equiv 0$. Therefore, since $\int_\Omega |\nabla u_n^+|^2\,dx \leq \beta \left[\int_\Omega (u_n^+)^q\,dx \right]^{\frac{2}{q}}$, we have $\lim_{n\rightarrow+\infty} \int_\Omega |\nabla u_n^+|^2\,dx = 0$, which leads to $\lim_{n\rightarrow+\infty} f_{\beta,\varepsilon}(u_n^+) = 0$, by applying the dominated convergence theorem. We then obtain a contradiction if we show that
$$
    \inf\{f_{\beta,\varepsilon}(w^+) ; w \in M^{\beta,\varepsilon}\} > 0 \,\,,\,\, \forall\varepsilon>0.
$$
Since for all $w\in M^{\beta,\varepsilon}$, $f_{\beta,\varepsilon}(w^+)=\max\{f_{\beta,\varepsilon}(tw^+) ; t>0\}$ (by Corollary \ref{cor:maxpositive}), it suffices to find two constants $\rho_{\beta,\varepsilon}>0$ and $c_{\beta,\varepsilon}>0$, such that for all $w\in H_0^1(\Omega)$ with $\int_\Omega |\nabla w|^2\,dx = \rho_{\beta,\varepsilon}$, we have $f_{\beta,\varepsilon}(w)\geq c_{\beta,\varepsilon}$. Indeed, if that last property is true, then for all $w\in M^{\beta,\varepsilon}$, since $w^+\not\equiv 0$, there exists $t_w>0$ such that $\int_\Omega |\nabla (t_w w^+)|^2\,dx = \rho_{\beta,\varepsilon}$, which leads to
$$
    f_{\beta,\varepsilon}(w^+)=\max\{f_{\beta,\varepsilon}(tw^+) ; t>0\} \geq f_{\beta,\varepsilon}(t_w w^+) \geq c_{\beta,\varepsilon}.
$$
Consider the sequence $(\mu_n)_n$ of (2,2)-eigenvalues.
For all $\phi\in H_0^1(\Omega)\setminus\{0\}$, H\"older's inequality gives
$$
    \frac{\norm{\nabla \phi}^2_{[L^2(\Omega)]^N}}{\norm{\phi}^2_{L^q(\Omega)}} \geq \frac{\norm{\nabla \phi}^2_{[L^2(\Omega)]^N}}{\norm{\phi}^2_{L^2(\Omega)}} |\Omega|^{\frac{q-2}{q}}\geq \mu_1|\Omega|^{\frac{q-2}{q}}
$$
Taking the infimum over all functions $\phi\in H_0^1(\Omega)\setminus\{0\}$, we obtain that $\lambda_1(\Omega)\geq \mu_1 |\Omega|^{1-\frac{2}{q}}$. Since $\beta>\lambda_1(\Omega)$ and $1<q<2$, we have $2^{2-\frac2q} \beta|\Omega|^{\frac2q-1}>\beta|\Omega|^{\frac2q-1} > \lambda_1(\Omega)|\Omega|^{\frac2q-1} \geq \mu_1$; hence, there exists $j_\beta\in\N$ such that $\mu_{j_\beta}\leq 2^{2-\frac{2}{q}}\beta|\Omega|^{\frac{2}{q}-1}<\mu_{j_\beta+1}$. We denote by $S_1$ and $S_2$ the subspaces of $H_0^1(\Omega)$ generated by the eigenfunctions of the Laplacian with Dirichlet boundary conditions associated to the eigenfunctions $\mu_j$, for $1\leq j\leq j_\beta$ and $j\geq j_\beta+1$ respectively.

Since $S_1$ is finite-dimensional and its elements are bounded (the Dirichlet eigenfunctions of the Laplacian belong to $L^\infty(\Omega)$, see \cite{Sun_Xing}), all norms are equivalent on $S_1$: there exists $C_\beta>0$ such that $\norm{w}_{L^\infty(\Omega)}\leq C_\beta\big(\int_\Omega|\nabla w|^2\,dx\big)^{\frac12}$ for all $w\in S_1$. Setting $\nu_{\beta,\varepsilon} = \frac{2^{\frac1q-1}\varepsilon}{C_\beta}$, we obtain $|w|\leq2^{\frac1q-1}\varepsilon$ a.e. in $\Omega$ for every $w\in S_1$ with $\int_\Omega|\nabla w|^2\,dx\leq\nu_{\beta,\varepsilon}^2$.

Let $w\in H_0^1(\Omega)$ such that $\int_\Omega |\nabla w|^2\,dx \leq \nu_{\beta,\varepsilon}^2$. Since the eigenfunctions of the Laplacian form a Hilbert basis of $H_0^1(\Omega)$, we can write $w=w_1+w_2$, with $w_1\in S_1$ and $w_2\in S_2$. As $w_1$ and $w_2$ belong to orthogonal subspaces of $H_0^1(\Omega)$), we have
\begin{equation}\label{Expression_f}
    f_{\beta,\varepsilon}(w) = f_{\beta,\varepsilon}(w_1+w_2) = \int_\Omega |\nabla w_1|^2\,dx + \int_\Omega |\nabla w_2|^2\,dx - q \left[\int_\Omega G_{\beta,\varepsilon}(w_1+w_2)\,dx\right]^\frac{2}{q}.
\end{equation}
Since $\int_\Omega |\nabla w_1|^2\,dx \leq \nu_{\beta,\varepsilon}^2$, then $|w_1|\leq 2^{\frac{1}{q}-1}\varepsilon$ a.e. in $\Omega$. For almost every $x\in\Omega$,  $|w_1+w_2|^q \leq 2^{q-1} (|w_1|^q+|w_2|^q)$. Therefore, on the set $\{w_1+w_2 > \varepsilon\} \subset\Omega$, we have
\[
\begin{split}
    2^{q-1}|w_2|^q &\geq |w_1+w_2|^q - 2^{q-1}|w_1|^q \geq |w_1+w_2|^q - \varepsilon^q \\
    &\geq |w_1+w_2|^q - \varepsilon^q - q\varepsilon^{q-1} |w_1+w_2| + q\varepsilon^q = q\int_\varepsilon^{|w_1+w_2|} (\tau^{q-1} - \varepsilon^{q-1})\,d\tau \\
    &\geq \beta^{-\frac{q}{2}} q^{\frac{q}{2}} \int_0^{|w_1+w_2|} g_{\beta,\varepsilon}(\tau) \,d\tau \, \chi_{\{w_1+w_2 \geq \varepsilon\}} = \beta^{-\frac{q}{2}} q^{\frac{q}{2}} G_{\beta,\varepsilon}(w_1+w_2).
\end{split}
\]
Moreover, this inequality is trivially satisfied on the set $\{w_1+w_2 \leq \varepsilon\}$ (because $G_{\beta,\varepsilon}(w_1+w_2) = 0$). Integrating this inequality on $\Omega$, we obtain
$$
    \int_\Omega G_{\beta,\varepsilon}(w_1+w_2)\,dx \leq 2^{q-1} \beta^{\frac{q}{2}} q^{-\frac{q}{2}} \int_\Omega |w_2|^q\,dx.
$$
Therefore,
\begin{equation}\label{Inegalite1}
    q \left[\int_\Omega G_{\beta,\varepsilon}(w_1+w_2)\,dx\right]^\frac{2}{q} \leq 2^{2-\frac{2}{q}} \beta \left[\int_\Omega |w_2|^q\,dx\right]^\frac{2}{q}\,.
\end{equation}
Since $w_2\in S_2$, by H\"older's inequality,
$$
  \int_\Omega |\nabla w_2|^2\,dx \geq \mu_{j_\beta + 1} \int_\Omega w_2^2\,dx \geq \mu_{j_\beta + 1} |\Omega|^{1-\frac{2}{q}} \left[\int_\Omega |w_2|^q\,dx\right]^{\frac{2}{q}},
$$
  and thus
$$ \left[\int_\Omega |w_2|^q\,dx\right]^{\frac{2}{q}} \leq \frac{|\Omega|^{\frac{2}{q}-1}}{\mu_{j_\beta + 1}} \int_\Omega |\nabla w_2|^2\,dx.
$$
Combining that last inequality with \eqref{Inegalite1}, we find
$$
    q \left[\int_\Omega G_{\beta,\varepsilon}(w_1+w_2)\,dx\right]^\frac{2}{q} \leq \frac{2^{2-\frac{2}{q}} \beta|\Omega|^{\frac{2}{q}-1}}{\mu_{j_\beta + 1}} \int_\Omega |\nabla w_2|^2\,dx.
$$
Consequently, using expression \eqref{Expression_f}, we obtain the following inequality:
$$
    f_{\beta,\varepsilon}(w) \geq \int_\Omega |\nabla w_1|^2\,dx + \left(1-\frac{2^{2-\frac{2}{q}}\beta|\Omega|^{\frac{2}{q}-1}}{\mu_{j_\beta+1}}\right)\int_\Omega |\nabla w_2|^2\,dx,
$$
where we notice that $\left(1-\frac{2^{2-\frac{2}{q}}\beta|\Omega|^{\frac{2}{q}-1}}{\mu_{j_\beta+1}}\right) > 0$, by the definition of $\mu_{j_\beta+1}$.

Hence, there exists $\tilde{c}_{\beta,\varepsilon}>0$ such that $f_{\beta,\varepsilon}(w) \geq \tilde{c}_{\beta,\varepsilon} \int_\Omega |\nabla w|^2\,dx$ for all $w\in H_0^1(\Omega)$ with $\int_\Omega |\nabla w|^2\,dx\leq \nu_{\beta,\varepsilon}^2$. Therefore, the constants $\rho_{\beta,\varepsilon}=\nu_{\beta,\varepsilon}^2$ and $c_{\beta,\varepsilon}=\tilde{c}_{\beta,\varepsilon}\,\nu_{\beta,\varepsilon}^2>0$ satisfy the desired properties. All of this implies $u^+ \not\equiv 0$.

By weak convergence in $H_0^1(\Omega)$, since $(u_n)_n \subset M^{\beta,\varepsilon}$, we have $f'_{\beta,\varepsilon}(u^+)[u^+] \leq \lim_{n\rightarrow +\infty} f'_{\beta,\varepsilon}(u_n^+)[u_n^+] = 0$. In addition, as was proved in Lemma \ref{IntermediateValueTheorem_f}, we have $\lim_{t\rightarrow 0} \frac{1}{t}f_{\beta,\varepsilon}'(tu^+)[u^+] = 2\int_\Omega |\nabla u^+|^2 \,dx > 0$ (because $u^+ \not\equiv 0$). Therefore, by continuity of $t\mapsto f'_{\beta,\varepsilon}(tu^+)[u^+]$, there exists $\tilde{t}\in(0,1]$ such that the function $\tilde{u} = -u^- + \tilde{t}u^+$ belongs to $M^{\beta,\varepsilon}$. Since for all $n\in\N$, $f_{\beta,\varepsilon}(\tilde{t}u_n^+) \leq f_{\beta,\varepsilon}(u_n^+)$, the weak $H_0^1$-convergence of $(-u_n^- + \tilde{t}u_n^+)_n$ to $\tilde{u}$ yields
$$
    f_{\beta,\varepsilon}(\tilde{u}) \leq \liminf_{n\rightarrow+\infty} f_{\beta,\varepsilon}(-u_n^- + \tilde{t}u_n^+) \leq \liminf_{n\rightarrow+\infty} f_{\beta,\varepsilon}(u_n) = \inf\{f_{\beta,\varepsilon}(u) ; u\in M^{\beta,\varepsilon}\}.
$$
Consequently, the infimum of $f_{\beta,\varepsilon}$ on $M^{\beta,\varepsilon}$ is achieved, and $f_{\beta,\varepsilon}(\tilde{u}) = \min_{M^{\beta,\varepsilon}} f_{\beta,\varepsilon}$.
\end{proof}

Having established the existence of a minimiser $u_{\beta,\varepsilon}$ of $f_{\beta,\varepsilon}$ on $M^{\beta,\varepsilon}$, the next step is to determine the Euler--Lagrange equation satisfied by $u_{\beta,\varepsilon}$, which is usually done using the Lagrange multipliers theorem. However, while the constraint $\int_\Omega(u^-)^q\,dx=1$ defines a $C^1$ manifold, the other constraint given by $f_{\beta,\varepsilon}'(u)[u^+] = 0$ does not provide enough regularity for the application of the theorem. Therefore, the strategy must be adapted to our case: in the following Proposition, we show directly that $u_{\beta,\varepsilon}$ is a critical point of $f_{\beta,\varepsilon}$ restricted to $S$.

\begin{prop}\label{PDEepsilon}
    Let $\beta>\lambda_1(\Omega)$, $\varepsilon>0$, and $u_{\beta,\varepsilon}$ be a minimiser of $f_{\beta,\varepsilon}$ on $M^{\beta,\varepsilon}$. Then $u_{\beta,\varepsilon}$ is a weak solution of the equation
    \begin{equation}
        \Delta u_{\beta,\varepsilon} + \left[\int_\Omega G_{\beta,\varepsilon}(u_{\beta,\varepsilon})\,dx\right]^{\frac{2}{q}-1} g_{\beta,\varepsilon}(u_{\beta,\varepsilon}) - \alpha_{\beta,\varepsilon} (u_{\beta,\varepsilon}^-)^{q-1} = 0 \,\,\,\,\textnormal{in } \Omega
    \end{equation}
    where $\alpha_{\beta,\varepsilon} = \int_\Omega |\nabla u_{\beta,\varepsilon}^-|^2\,dx$.
\end{prop}

\begin{proof}
We aim to show that there exists a Lagrange multiplier $\alpha_{\beta,\varepsilon}\in\R$ such that
\begin{equation}\label{EquationEps}
    \frac{1}{2}f'_{\beta,\varepsilon}(u_{\beta,\varepsilon})[\psi] = -\alpha_{\beta,\varepsilon} \int_\Omega (u_{\beta,\varepsilon}^-)^{q-1} \psi\,dx \,\,,\,\, \forall\psi\in H_0^1(\Omega).
\end{equation}
By Corollary \ref{cor:maxpositive}, $u_{\beta,\varepsilon}$ is the unique maximum point of $f_{\beta,\varepsilon}$ on the set $\{u_{\beta,\varepsilon} + tu_{\beta,\varepsilon}^+\,;\, t\geq -1\}$. Note that the constraint defining $S$, that is, $J(u)=\int_\Omega(u^-)^q\,dx$, is of class $C^1$ on $H_0^1(\Omega)$ for $q>1$, and $J'(u) \neq 0$ on $S$, so that $S$ is a $C^1$ submanifold of $H_0^1(\Omega)$, while $f_{\beta,\varepsilon}\in C^1(H_0^1(\Omega))$.

Suppose that \eqref{EquationEps} fails for every multiplier. Set $u=u_{\beta,\varepsilon}$, $f=f_{\beta,\varepsilon}$ and $c=f(u)$. On the open set $\{z:\norm{z^-}_{L^q(\Omega)}>0\}$, define the $C^1$ normalisation
\[
    R \,:\, z \, \mapsto \, \frac{z}{\norm{z^-}_{L^q(\Omega)}}.
\]
It maps into $S$ and equals the identity on $S$. Since $f|_S$ is not critical at $u$, there exists $v\in H_0^1(\Omega)$ such that the Gateaux derivative $(f\circ R)'(u)[v]<0$. Set $p(t)=-u^-+(1+t)u^+$ for $|t|\leq\frac12$. By continuity, we can choose $d\in(0,\frac12)$, $h>0$ and $a>0$ such that
\[
    \frac{\partial}{\partial s}f\big(R(p(t)+sv)\big)\leq-a \quad\text{for }|t|\leq d,\quad 0\leq s\leq h.
\]
Shrinking $d$ and $h$ also ensures that all these functions have non-zero positive parts. Choose a continuous function $\rho:[-\frac12,\frac12]\to[0,1]$ with $\rho(0)=1$ and $\rho(t)=0$ for $|t|\geq d$, and set $\eta(t)=R(p(t)+h\rho(t)v)$. Then $\eta$ is continuous, belongs to $S$, has a non-zero positive part and agrees with $p$ at the endpoints. Moreover,
\[
    f(\eta(t))\leq f(p(t))-ah\rho(t)<c \quad\text{for }|t|\leq\tfrac{1}{2},
\]
where strictness at $t=0$ follows from $\rho(0)=1$, and strictness for $t\neq0$ from Corollary \ref{cor:maxpositive}. This gives the required local deformation directly; for general deformation methods on $C^1$ manifolds, see \cite{Willem,Szulkin}.
Consider the function $\varphi(t)=f_{\beta,\varepsilon}'(\eta(t))[\eta(t)^+]$,
which is continuous on $[-\tfrac12,\tfrac12]$.
Since $\eta(\pm\tfrac12)=u_{\beta,\varepsilon}\pm\tfrac12 u_{\beta,\varepsilon}^+
=-u_{\beta,\varepsilon}^-+(1\pm\tfrac12)u_{\beta,\varepsilon}^+$, Lemma \ref{Lemma1} yields $\varphi(-\tfrac12)>0$ and $\varphi(\tfrac12)<0$.
By the intermediate value theorem, there exists $\bar{t}\in[-\tfrac{1}{2},\tfrac{1}{2}]$ such that $\varphi(\bar t)=0$; since $\norm{\eta(\bar t)^-}_{L^q(\Omega)}=1$ and $\eta(\bar t)^+\not\equiv0$, this means $\eta(\bar{t}) \in M^{\beta,\varepsilon}$. This yields a contradiction, because $f_{\beta,\varepsilon}(\eta(\bar{t})) < f_{\beta,\varepsilon}(u_{\beta,\varepsilon}) = \min_{M^{\beta,\varepsilon}} f_{\beta,\varepsilon}$. Therefore, there exists a Lagrange multiplier $\alpha_{\beta,\varepsilon}\in\R$ such that equation \eqref{EquationEps} is satisfied.

By taking $\psi=u_{\beta,\varepsilon}^- \in H_0^1(\Omega)$, we obtain $\alpha_{\beta,\varepsilon} = -\frac{1}{2}f'_{\beta,\varepsilon}(u_{\beta,\varepsilon})[u_{\beta,\varepsilon}^-] = \int_\Omega |\nabla u_{\beta,\varepsilon}^-|^2\,dx$.
\end{proof}

\subsection{Taking the limit $\varepsilon\rightarrow 0$ for the family $(f_{\beta,\varepsilon})_\varepsilon$}

To summarise the previous subsection, we have minimised the functional $f_{\beta,\varepsilon}$ for $\varepsilon>0$ and $\beta>\lambda_1(\Omega)$, and found the equation satisfied by the minimiser $u_{\beta,\varepsilon}$. Let us remark that we have studied the regularised problem with $\varepsilon$ instead of directly considering the limit case $\varepsilon=0$ because for $\varepsilon>0$ the non-linearity $g_{\beta,\varepsilon}$ vanishes on $[0,\varepsilon]$, so that $f_{\beta,\varepsilon}(tw)>0$ for $t>0$ small (Lemma \ref{IntermediateValueTheorem_f}) and $\inf\{f_{\beta,\varepsilon}(w^+)\,;\,w\in M^{\beta,\varepsilon}\}>0$ (proof of Proposition \ref{InfAchieved})
We must now let $\varepsilon$ go to $0$ to obtain the desired minimality property and equation satisfied by the limit function.

We start this subsection by an important Lemma that allows us to decrease the gradients of two functions supported on disjoint subsets of $\Omega$, by a slight modification of the supports of the functions.

\begin{Lemma}\label{ModificationOfSupports}
    Let $w\in H_0^1(\Omega)$. Suppose that $\norm{w^+}_{L^q(\Omega)} = \norm{w^-}_{L^q(\Omega)} = 1$. Moreover, suppose there exists two constants $A,B>0$ such that
    $$
        \int_\Omega |\nabla w^+|^2 \,dx \leq A \,\, \textnormal{, and } \,\, \int_\Omega |\nabla w^-|^2 \,dx \leq B.
    $$
    Then, if one of those inequalities is strict, there exists $v\in H_0^1(\Omega)$ such that $\norm{v^+}_{L^q(\Omega)} = \norm{v^-}_{L^q(\Omega)} = 1$,
    $$
        \int_\Omega |\nabla v^+|^2 \,dx < A \,\, \textnormal{, and } \,\, \int_\Omega |\nabla v^-|^2 \,dx < B.
    $$
\end{Lemma}

\begin{proof}
    If both inequalities are strict, then the function $v=w$ is suitable. If not, without loss of generality, we suppose that only the second inequality is strict. For $t > 0$, we set
    $$
        \phi_t =\frac{(w^- - t)_+}{\norm{(w^- - t)_+}_{L^q(\Omega)}}.
    $$
    The sequence $(\phi_t)_t$ converges to $w^-$ in $H_0^1(\Omega)$, as $t\to 0$. Hence, since $\int_\Omega |\nabla w^-|^2 \,dx < B$, there exists $t_0>0$ such that for all $t<t_0$, we have $\int_\Omega |\nabla \phi_t|^2 \,dx < B$.

Let $t\in(0,t_0)$. There exist $r,s$ such that $0<r<s<t$ and $|\{r<w^-<s\}|>0$. Otherwise $|\{0<w^-<t\}|=0$, so that $\min\{w^-,t\}=t\,\chi_{\{w^-\geq t\}}\in H_0^1(\Omega)$ has zero gradient a.e. in $\Omega$ (since $\nabla w^-=0$ a.e. on $\{w^-=0\}$) and is non-zero (since $w^-\not\equiv0$), which contradicts the Poincar\'e inequality.
We then define the function $\zeta_t = \min\{(w^- - r)_+, (s - w^-)_+\}$. Then $\zeta_t \in H_0^1(\Omega)$, $\zeta_t \geq 0$, $\zeta_t \not \equiv 0$, and its support is contained in $\{r < w^- < s\}$. Consequently, for all $t\in(0,t_0)$, we have
$$
    w^+ \zeta_t = 0, \qquad \phi_t \zeta_t = 0, \qquad w^+ \phi_t =0, \quad\textnormal{a.e. in }\Omega.
$$
Fix $t\in(0,t_0)$ and the corresponding $r,s,\zeta_t$ above. For an independent parameter $\delta>0$, set
\[
    \psi_\delta=\frac{w^++\delta\zeta_t}{\norm{w^++\delta\zeta_t}_{L^q(\Omega)}}.
\]
Write $a=\int_\Omega|\nabla w^+|^2\,dx>0$, $b=\int_\Omega|\nabla\zeta_t|^2\,dx$ and $c=\norm{\zeta_t}_{L^q(\Omega)}^q>0$. Since $w^+\zeta_t=0$ a.e.,
\[
    \int_\Omega|\nabla\psi_\delta|^2\,dx = \frac{a+\delta^2b}{(1+\delta^q c)^{\frac{2}{q}}}.
\]
As $q<2$, choose $\delta>0$ so small that $\delta^{2-q}b < \frac{2ac}{q}$. Bernoulli's inequality gives $(1+\delta^q c)^{\frac{2}{q}}\geq 1 + \frac{2}{q}\delta^q c$, hence $\int_\Omega|\nabla\psi_\delta|^2\,dx<a\leq A$. Finally, $\psi_\delta\phi_t=0$ a.e.; therefore $v=\psi_\delta-\phi_t$ satisfies $\norm{v^+}_{L^q(\Omega)}=\norm{v^-}_{L^q(\Omega)}=1$, $\int_\Omega|\nabla v^+|^2\,dx<A$ and $\int_\Omega|\nabla v^-|^2\,dx<B$.
\end{proof}

\begin{prop}\label{ConvergenceSequence}
    For all $\beta>\lambda_1(\Omega)$ and $\varepsilon>0$, let $u_{\beta,\varepsilon}$ be a minimiser of $f_{\beta,\varepsilon}$ on $M^{\beta,\varepsilon}$. Let $\bar{u}_{\beta,\varepsilon} = -u_{\beta,\varepsilon}^- + \frac{u_{\beta,\varepsilon}^+}{\norm{u_{\beta,\varepsilon}^+}_{L^q(\Omega)}}$. Therefore, up to  a subsequence,  $(\bar{u}_{\beta,\varepsilon})_\varepsilon$ converges in $H_0^1(\Omega)$, as $\varepsilon$ tends to $0$, to a function $\bar{u}_\beta$ such that $\int_\Omega |\nabla \bar{u}_\beta^+|^2\,dx = \beta$ and
    \[
\begin{split}\int_\Omega |\nabla \bar{u}_\beta^-|^2\,dx = \min\Big\{\int_\Omega |\nabla u^-|^2\,dx \,;\ & u\in H_0^1(\Omega),\ \norm{u^+}_{L^q(\Omega)} = \norm{u^-}_{L^q(\Omega)} = 1 ,\\ & \int_\Omega |\nabla u^+|^2\,dx=\beta \Big\}.\end{split}
\]
\end{prop}

\begin{proof}
Let $\beta>\lambda_1(\Omega)$. Since $u_{\beta,\varepsilon}\in M^{\beta,\varepsilon}$, we have for all $\varepsilon>0$,
$$
    \int_\Omega|\nabla u_{\beta,\varepsilon}^+|^2\,dx = \left[\int_\Omega G_{\beta,\varepsilon}(u_{\beta,\varepsilon})\,dx\right]^{\frac{2}{q}-1}\int_\Omega g_{\beta,\varepsilon}(u_{\beta,\varepsilon})u_{\beta,\varepsilon}^+\,dx < \beta \left[\int_\Omega(u_{\beta,\varepsilon}^+)^q\,dx\right]^{\frac{2}{q}}.
$$
Then, for all $\varepsilon>0$, $\int_\Omega |\nabla \bar{u}_{\beta,\varepsilon}^+|^2\,dx < \beta$.

On the other hand, $\limsup_{\varepsilon\rightarrow0} \int_\Omega |\nabla u _{\beta,\varepsilon}^-|^2\,dx<+\infty$. Indeed, since $\beta>\lambda_1(\Omega)$, there exists $\hat{u}\in H_0^1(\Omega)$ such that $\int_\Omega (\hat{u}^+)^q\,dx = \int_\Omega (\hat{u}^-)^q \,dx = 1$, and $\int_\Omega |\nabla \hat{u}^+|^2\,dx < \beta$ (see Remark \ref{rem:testfunction}).
Since $u_{\beta,\varepsilon}$ minimises $f_{\beta,\varepsilon}$ on $M^{\beta,\varepsilon}$, for all $\varepsilon>0$,
$$
    f_{\beta,\varepsilon}(u_{\beta,\varepsilon}) \leq \int_\Omega |\nabla \hat{u}^-|^2\,dx + \max\{f_{\beta,\varepsilon}(t\hat{u}^+) ; t\geq 0\},
$$
because the real number $t_\varepsilon$ for which the maximum is attained verifies $-\hat{u}^-+t_\varepsilon \hat{u}^+ \in M^{\beta,\varepsilon}$.  By Lemma \ref{IntermediateValueTheorem_f}, we have $\lim_{\varepsilon\rightarrow 0} \max\{f_{\beta,\varepsilon}(t\hat{u}^+) ; t\geq 0\} =0$.  Since $u_{\beta,\varepsilon}\in M^{\beta,\varepsilon}$, by Corollary \ref{cor:maxpositive},
$$
    \limsup_{\varepsilon\rightarrow0} \int_\Omega |\nabla u _{\beta,\varepsilon}^-|^2\,dx   \leq \limsup_{\varepsilon\rightarrow0} f_{\beta,\varepsilon}(u_{\beta,\varepsilon}) \leq \int_\Omega |\nabla \hat{u}^-|^2\,dx < +\infty.
$$
Up to extraction, the sequence $(\bar{u}_{\beta,\varepsilon})_\varepsilon$ converges as $\varepsilon$ tends to $0$ to a function $\bar{u}_\beta \in H_0^1(\Omega)$, strongly in $L^q(\Omega)$, weakly in $H_0^1(\Omega)$, and a.e. in $\Omega$.

We want to show that the sequence $(\bar{u}_{\beta,\varepsilon})_\varepsilon$ actually converges  in $H_0^1(\Omega)$. It suffices to show the convergence of the norms
\begin{equation}\label{SupportsEquality}
    \int_\Omega |\nabla \bar{u}_\beta^+|^2 = \beta \textnormal{ , and } \int_\Omega |\nabla \bar{u}_\beta^-|^2\,dx = \lim_{\varepsilon\rightarrow0} \int_\Omega |\nabla u_{\beta,\varepsilon}^-|^2\,dx.
\end{equation}
By contradiction, we assume that $\int_\Omega |\nabla \bar{u}_\beta^+|^2 < \beta$, or $\int_\Omega |\nabla \bar{u}_\beta^-|^2\,dx < \lim_{\varepsilon\rightarrow0} \int_\Omega |\nabla u_{\beta,\varepsilon}^-|^2\,dx$ (up to extraction).  By Lemma \ref{ModificationOfSupports}, we can construct a function $\tilde{u}_\beta \in H_0^1(\Omega)$ such that $\norm{\tilde{u}_\beta^+}_{L^q(\Omega)} = \norm{\tilde{u}_\beta^-}_{L^q(\Omega)} = 1$, $\int_\Omega |\nabla \tilde{u}_\beta^+|^2\,dx < \beta$, and $\int_\Omega |\nabla \tilde{u}_\beta^-|^2\,dx < \lim_{\varepsilon\rightarrow0} \int_\Omega |\nabla u_{\beta,\varepsilon}^-|^2\,dx$.

For all $\varepsilon>0$, consider $\tilde{u}_{\beta,\varepsilon} \in M^{\beta,\varepsilon}$ such that $\tilde{u}_{\beta,\varepsilon}^- = \tilde{u}_\beta^-$, and $\tilde{u}_{\beta,\varepsilon}^+ = t_\varepsilon \tilde{u}_\beta^+$ where $t_\varepsilon$ is the unique positive real number such that $f_{\beta,\varepsilon}'(t_\varepsilon \tilde{u}_\beta^+)[\tilde{u}_\beta^+] = 0$ (such a value exists because $\int_\Omega |\nabla \tilde{u}_\beta^+|^2\,dx<\beta$, by Lemma \ref{IntermediateValueTheorem_f}). Therefore,
$$
    f_{\beta,\varepsilon}(u_{\beta,\varepsilon}) - f_{\beta,\varepsilon}(\tilde{u}_{\beta,\varepsilon}) = f_{\beta,\varepsilon}(u_{\beta,\varepsilon}^+) - f_{\beta,\varepsilon}(\tilde{u}_{\beta,\varepsilon}^+) + f_{\beta,\varepsilon}(-u_{\beta,\varepsilon}^-) - f_{\beta,\varepsilon}(-\tilde{u}_{\beta,\varepsilon}^-),
$$
where for all $\varepsilon>0$,
$f_{\beta,\varepsilon}(u_{\beta,\varepsilon}^+) \geq 0$,
$\lim_{\varepsilon\rightarrow0} f_{\beta,\varepsilon}(\tilde{u}_{\beta,\varepsilon}^+) \leq 0$
(since $\int_\Omega |\nabla \tilde{u}_\beta^+|^2\,dx < \beta$),
and
$\lim_{\varepsilon\rightarrow0} f_{\beta,\varepsilon}(-u_{\beta,\varepsilon}^-) = \lim_{\varepsilon\rightarrow0} \int_\Omega |\nabla u_{\beta,\varepsilon}^-|^2\,dx > \int_\Omega |\nabla \tilde{u}_{\beta}^-|^2\,dx = f_{\beta,\varepsilon}(-\tilde{u}_{\beta,\varepsilon}^-)$.
Consequently, for $\varepsilon>0$ sufficiently small, $f_{\beta,\varepsilon}(u_{\beta,\varepsilon}) > f_{\beta,\varepsilon}(\tilde{u}_{\beta,\varepsilon})$.
This contradicts the fact that $f_{\beta,\varepsilon}(u_{\beta,\varepsilon}) = \min_{M^{\beta,\varepsilon}} f_{\beta,\varepsilon}$, and $\tilde{u}_{\beta,\varepsilon} \in M^{\beta,\varepsilon}$.
This proves the convergence in $H_0^1(\Omega)$ of the sequence $(\bar{u}_{\beta,\varepsilon})_\varepsilon$ to $\bar{u}_\beta$, and $\int_\Omega |\nabla \bar{u}^+_\beta|^2\,dx = \beta$.

It remains to show that
\[
\begin{split}\int_\Omega |\nabla \bar{u}_\beta^-|^2\,dx = \min\Big\{\int_\Omega |\nabla u^-|^2\,dx \,;\ & u\in H_0^1(\Omega),\ \norm{u^+}_{L^q(\Omega)} = \norm{u^-}_{L^q(\Omega)} = 1 ,\\ & \int_\Omega |\nabla u^+|^2\,dx=\beta \Big\}.\end{split}
\]
If $v\in H_0^1(\Omega)$ satisfies $\norm{v^+}_{L^q(\Omega)}=\norm{v^-}_{L^q(\Omega)}=1$, $\int_\Omega|\nabla v^+|^2\,dx=\beta$ and $\int_\Omega|\nabla v^-|^2\,dx<\int_\Omega|\nabla\bar u_\beta^-|^2\,dx$, Lemmas \ref{ModificationOfSupports} and \ref{IntermediateValueTheorem_f} provide, exactly as above, functions $\hat v_{\beta,\varepsilon}\in M^{\beta,\varepsilon}$ with $f_{\beta,\varepsilon}(\hat v_{\beta,\varepsilon})<f_{\beta,\varepsilon}(u_{\beta,\varepsilon})$ for $\varepsilon$ small, a contradiction.
\end{proof}

\begin{prop}\label{PDE_limit}
    For all $\beta>\lambda_1(\Omega)$ and $\varepsilon>0$, let $u_{\beta,\varepsilon}$ be a minimiser of $f_{\beta,\varepsilon}$ on $M^{\beta,\varepsilon}$. Therefore, up to  a subsequence,  $(u_{\beta,\varepsilon})_\varepsilon$ converges in $H_0^1(\Omega)$, as $\varepsilon \to 0$, to a sign-changing function $u_\beta$ satisfying
    $$
        \Delta u_\beta - \alpha_\beta (u_\beta^-)^{q-1} + \beta \norm{u_\beta^+}_{L^q(\Omega)}^{2-q} (u_\beta^+)^{q-1} = 0 \,\,\,\,\,\,\ \textnormal{ in } \Omega,
    $$
    where
    \begin{equation}\label{alpha_beta}
        \alpha_\beta = \min\left\{\int_\Omega |\nabla u^-|^2\,dx \,;\, u\in H_0^1(\Omega), \norm{u^+}_{L^q(\Omega)} = \norm{u^-}_{L^q(\Omega)} = 1 , \int_\Omega |\nabla u^+|^2\,dx=\beta \right\}.
    \end{equation}
    Moreover, the function $\bar{u}_{\beta} = -u_{\beta}^- + \frac{u_{\beta}^+}{\norm{u_{\beta}^+}_{L^q(\Omega)}}$ is a minimiser for $\alpha_\beta$, and $(\alpha_\beta,\beta)\in\Sigma$.
\end{prop}

\begin{proof}
By Proposition \ref{PDEepsilon}, for all $\beta>\lambda_1(\Omega)$ and $\varepsilon>0$, $u_{\beta,\varepsilon}$ is a weak solution of
\begin{equation}\label{PDE_epsilon}
    \Delta u_{\beta,\varepsilon} + \left[\int_\Omega G_{\beta,\varepsilon}(u_{\beta,\varepsilon})\,dx\right]^{\frac{2}{q}-1} g_{\beta,\varepsilon}(u_{\beta,\varepsilon}) - \alpha_{\beta,\varepsilon} (u_{\beta,\varepsilon}^-)^{q-1} = 0 \,\,\,\,\textnormal{in } \Omega,
\end{equation}
where $\alpha_{\beta,\varepsilon} = \int_\Omega |\nabla u_{\beta,\varepsilon}^-|^2\,dx$. Moreover, by Proposition \ref{ConvergenceSequence}, the sequence $(\bar{u}_{\beta,\varepsilon})_\varepsilon$ converges in $H_0^1(\Omega)$ as $\varepsilon$ tends to $0$ to a function $\bar{u}_\beta$.

We show first that for all $\beta>\lambda_1(\Omega)$, $\liminf_{\varepsilon\rightarrow0} \norm{u_{\beta,\varepsilon}^+}_{L^q(\Omega)} > 0$. By contradiction, we assume that there exists $\beta>\lambda_1(\Omega)$ such that, up to extraction, $\lim_{\varepsilon\rightarrow0} \norm{u_{\beta,\varepsilon}^+}_{L^q(\Omega)} = 0$. Therefore, the sequence $(u_{\beta,\varepsilon}^+)_\varepsilon$ converges a.e. to $0$. In addition, since $f_{\beta,\varepsilon}'(u_{\beta,\varepsilon})[u_{\beta,\varepsilon}^+] = 0$, we obtain that $\lim_{\varepsilon\rightarrow0} \int_\Omega |\nabla u_{\beta,\varepsilon}^+|^2\,dx = 0$. Consequently, $(u_{\beta,\varepsilon}^+)_\varepsilon$ converges in $H_0^1(\Omega)$ to $0$.  Hence $\bar u_\beta^-$ is a non-negative $(2,q)$-eigenfunction on $\Omega$, and
Proposition \ref{FirstEigenFunctionsArePositive} yields $\bar{u}_\beta^- > 0$ a.e. in $\Omega$. On the other hand, by Proposition \ref{ConvergenceSequence}, $\norm{\bar{u}_\beta^+}_{L^q(\Omega)}=1$ and $\bar{u}_\beta^+ \bar{u}_\beta^- = 0$ a.e.\ in $\Omega$, so $\bar{u}_\beta^-$ vanishes on $\{\bar{u}_\beta^+ > 0\}$, which has positive measure, leading to a contradiction.

We then show that for all $\beta>\lambda_1(\Omega)$, $\limsup_{\varepsilon\rightarrow0} \norm{u_{\beta,\varepsilon}^+}_{L^q(\Omega)} < +\infty$.
By contradiction, we assume that there exists $\beta>\lambda_1(\Omega)$
such that, up to a subsequence, $\lim_{\varepsilon\rightarrow0} \norm{u_{\beta,\varepsilon}^+}_{L^q(\Omega)} = +\infty$.
By dividing equation \eqref{PDE_epsilon} by $t_\varepsilon = \norm{u_{\beta,\varepsilon}^+}_{L^q(\Omega)}$ and letting $\varepsilon$ go to $0$
(notice that $\alpha_{\beta,\varepsilon}$ is bounded, and $t_\varepsilon \rightarrow +\infty$),
we obtain $\int_\Omega[\nabla\bar{u}_\beta^+\cdot\nabla\psi - \beta(\bar{u}_\beta^+)^{q-1}\psi] \,dx = 0$ for all $\psi\in H_0^1(\Omega)$.
Since $\bar{u}_\beta^+ \geq 0$ and $\norm{\bar{u}_\beta^+}_{L^q(\Omega)} = 1$, by Proposition \ref{PositiveEigenfunctionsAreFirst}, $\beta = \lambda_1(\Omega)$, which is a contradiction.

After extracting a further subsequence, $\norm{u_{\beta,\varepsilon}^+}_{L^q(\Omega)}$ converges to a strictly positive real number.
Consequently, for all $\beta>\lambda_1(\Omega)$, the sequence $(u_{\beta,\varepsilon})_\varepsilon$ converges in $H_0^1(\Omega)$ to a sign-changing function $u_\beta$. By letting $\varepsilon$ go to $0$ in equation \eqref{PDE_epsilon}, we observe that $u_\beta$ is a weak solution of
$$
    \left\{\begin{array}{cl}
        \Delta u_\beta - \alpha_\beta (u_\beta^-)^{q-1} + \beta \norm{u_\beta^+}_{L^q(\Omega)}^{2-q} (u_\beta^+)^{q-1} = 0 & \,\,\,\, \textnormal{on } \Omega \\
        u_\beta = 0 & \,\,\,\, \textnormal{in } \partial\Omega.
    \end{array}
    \right.
$$
Thus, $(\alpha_\beta,\beta)\in\Sigma$. Moreover, we have
\[
\begin{split}
    \alpha_\beta &= \lim_{\varepsilon\rightarrow0} \alpha_{\beta,\varepsilon} = \lim_{\varepsilon\rightarrow0} \int_\Omega |\nabla \bar{u}_{\beta,\varepsilon}^-|^2\,dx = \int_\Omega |\nabla \bar{u}_{\beta}^-|^2\,dx = \int_\Omega |\nabla u_{\beta}^-|^2\,dx \\
    & = \min\left\{\int_\Omega |\nabla u^-|^2\,dx \,;\, u\in H_0^1(\Omega), \norm{u^+}_{L^q(\Omega)} = \norm{u^-}_{L^q(\Omega)} = 1 , \int_\Omega |\nabla u^+|^2\,dx=\beta \right\}.
\end{split}
\]
Therefore, for all $\beta>\lambda_1(\Omega)$,
$$
\int_\Omega |\nabla u_\beta
^+|^2\,dx = \beta \left[\int_\Omega (u_\beta^+)^{q}\,dx\right]^{\frac{2}{q}} > 0\,,\quad
\int_\Omega |\nabla u_\beta^-|^2\,dx = \alpha_\beta\,.
$$
Since $u_\beta^+ \not\equiv 0$, then $\alpha_\beta = \int_\Omega |\nabla u_\beta^-|^2\,dx > \lambda_1(\Omega)$; indeed, if $\int_\Omega |\nabla u_\beta^-|^2\,dx = \lambda_1(\Omega) = \lambda_1(\Omega) \norm{u_\beta^-}_{L^q(\Omega)}^{2}$, $u_\beta^-$ would be a $(2,q)$-eigenfunction associated with $\lambda_1(\Omega)$, and thus, by Proposition \ref{FirstEigenFunctionsArePositive}, we would have $u_\beta^- > 0$ a.e. in $\Omega$, contradicting $u_\beta^+ \not\equiv 0$.
\end{proof}

By combining all of the above propositions, we have obtained the following theorem.

\begin{thm}\label{FucikSelectedEigenfunction}
    Let $\Omega$ be a bounded connected open set in $\R^N$. For all $\beta>\lambda_1(\Omega)$, set
    \begin{equation}\label{Definition_alpha_beta}
        \alpha_\beta = \inf\left\{\int_\Omega |\nabla u^-|^2\,dx \,;\, u\in H_0^1(\Omega), \norm{u^+}_{L^q(\Omega)} = \norm{u^-}_{L^q(\Omega)} = 1 , \int_\Omega |\nabla u^+|^2\,dx=\beta \right\}.
    \end{equation}
   For all $\beta>\lambda_1(\Omega)$, we have $\alpha_\beta > \lambda_1(\Omega)$, and $(\alpha_\beta,\beta)\in\Sigma$. Moreover, the infimum in the definition of $\alpha_\beta$ is achieved by the function
   $\bar{u}_{\beta} = -u_{\beta}^- + \frac{u_{\beta}^+}{\norm{u_{\beta}^+}_{L^q(\Omega)}}$, where $u_\beta$ solves the Fu\v{c}\'\i k eigenvalue problem
    \begin{equation}\label{Fucik_alpha_beta}
        \left\{
        \begin{array}{cl}
            -\Delta u_\beta = - \alpha_\beta (u_\beta^-)^{q-1} + \beta\norm{u_\beta^+}_{L^q(\Omega)}^{2-q}(u_\beta^+)^{q-1} & \textnormal{ in } \Omega \\
            u_\beta = 0 & \textnormal{ on } \partial\Omega.
        \end{array}
        \right.
    \end{equation}
\end{thm}

\begin{defn}
    Let $\Omega$ be a bounded connected open set in $\R^N$. For all $\beta > \lambda_1(\Omega)$, we shall call the function $u_\beta$ defined in Theorem \ref{FucikSelectedEigenfunction}, the \textbf{Fu\v{c}\'\i k-selected $(2,q)$-eigenfunction} associated to the pair $(\alpha_\beta,\beta)\in\Sigma$.
\end{defn}

\begin{rem}
    As stated in Remark \ref{RescalingEigenfunctions}, we may not assume that $\norm{u_\beta^+}_{L^q(\Omega)} = \norm{u_\beta^-}_{L^q(\Omega)} = 1$. Therefore, we cannot conclude that the minimiser $\bar{u}_\beta$ of the definition \eqref{Definition_alpha_beta} of $\alpha_\beta$ yields a solution of the following equation
    $$
        -\Delta v = - \alpha_\beta (v^-)^{q-1} + \beta(v^+)^{q-1}.
    $$
\end{rem}

\subsection{Properties of $(\alpha_\beta,\beta)$}

By defining for all $\beta>\lambda_1(\Omega)$, a real value $\alpha_\beta > \lambda_1(\Omega)$ such that $(\alpha_\beta,\beta)\in\Sigma$, we have constructed a curve included in the Fu\v{c}\'\i k $(2,q)$-spectrum. In other words, the application
\[
\begin{split}
    \begin{array}{rccc}
        A : &(\lambda_1(\Omega),+\infty) & \rightarrow & (\lambda_1(\Omega),+\infty) \\
        &\beta & \mapsto & \alpha_\beta
    \end{array}
\end{split}
\]
is well-defined. In this section we prove some properties of this application.
We start by proving the following Proposition, which explains why the associated curve is called \textbf{the first non-trivial curve} of the Fu\v{c}\'\i k spectrum.

\begin{prop}\label{SmallestFucikPair}
    For all $\beta>\lambda_1(\Omega)$, and all $\alpha>\lambda_1(\Omega)$ such that $(\alpha,\beta)\in\Sigma$, we have $\lambda_1(\Omega)<\alpha_\beta\leq\alpha$.
\end{prop}

\begin{proof}
Let $\beta>\lambda_1(\Omega)$, and let $\alpha>\lambda_1(\Omega)$ such that $(\alpha,\beta) \in \Sigma$.  Then, there exists $u\in H_0^1(\Omega)$ that is a weak solution of
$$
    \left\{\begin{array}{cl}
        -\Delta u = -\alpha \norm{u^-}_{L^q(\Omega)}^{2-q} (u^-)^{q-1} + \beta \norm{u^+}_{L^q(\Omega)}^{2-q} (u^+)^{q-1} & \,\,\,\, \textnormal{in } \Omega \\
        u=0 & \,\,\,\, \textnormal{on } \partial\Omega.
    \end{array}
    \right.
$$
Multiplying the equation by $u^+$ and $u^-$ and integrating on $\Omega$ yields
$$
    \int_\Omega |\nabla u^+|^2 \,dx = \beta \norm{u^+}_{L^q(\Omega)}^2 \,, \,\,\,\,\, \int_\Omega |\nabla u^-|^2 \,dx = \alpha \norm{u^-}_{L^q(\Omega)}^2.
$$
We consider the nodally normalised function
$$
    v = \frac{u^+}{\norm{u^+}_{L^q(\Omega)}} - \frac{u^-}{\norm{u^-}_{L^q(\Omega)}}.
$$
Although $v$ might not be a solution to the equation, we still have $\norm{v^+}_{L^q(\Omega)} = \norm{{v^-}}_{L^q(\Omega)} = 1$, $\int_\Omega |\nabla v^+|^2 \,dx = \beta$, and $\int_\Omega |\nabla {v^-}|^2 \,dx = \alpha$. Therefore, by the variational definition of $\alpha_\beta$, we have $\alpha_\beta \leq \int_\Omega |\nabla {v^-}|^2\,dx = \alpha$.
\end{proof}

The next lemma is used in the proof of the strict monotonicity of $\beta\mapsto\alpha_\beta$ (Proposition \ref{CurveIsContinuousAndDecreasing}) and in Proposition \ref{Gamma = gamma}. In the case $q=2$ the analogous statement is used in \cite{Molle_Passaseo}; we include a short proof based on periodic homogenisation.

\begin{Lemma}\label{DecreasingTheNorm}
    Let $\beta>0$. Let $w\in H_0^1(\Omega)$ such that $w\geq 0$, $w\not\equiv0$, and $\int_\Omega|\nabla w|^2\,dx \leq \beta\|w\|_{L^q(\Omega)}^2$. Then, there exists $v\in H_0^1(\Omega)\setminus\{0\}$ such that $0\leq {v} \leq w$ and $\int_\Omega |\nabla v|^2\,dx = \beta \|v\|_{L^q(\Omega)}^2$.
\end{Lemma}

\begin{proof}
    If equality holds, take $v=w$. Otherwise, fix a $\mathbb{Z}^N$-periodic $\chi\in C^\infty(\R^N)$, such that $0\le\chi\le1$ and that $\chi$ is non-constant. For $k\in\N$ and $t\in[0,1]$, set $v_{t,k}=w\,[(1-t)+t\chi(kx)]$, so $0\leq v_{t,k}\le w$, and $v_{t,k}\not\equiv0$. The Rayleigh quotient $Q/(t,k) \mapsto \frac{\int_\Omega|\nabla v_{t,k}|^2\,dx}{\norm{v_{t,k}}_{L^q(\Omega)}^2}$ is continuous in $t$, and verifies $Q(0,k)<\beta$.

    A homogenisation lemma shows that, if $h\in L^1(\Omega)$ and if $F\in L^\infty((0,1)^N)$ is $\Z^N$-periodic, then
    $$
        \int_\Omega h(x)F(kx) \,dx \xrightarrow[k\rightarrow+\infty]{} \left[ \int_{(0,1)^N} F(\tau) \,d\tau \right] \left[ \int_\Omega h(x) \,dx \right].
    $$
    Thus, $\int|\nabla(w\chi(k\cdot))|^2\,dx\sim k^2\int w^2\,dx\cdot\bar{c}_\chi\to+\infty$ (where $\bar{c}_\chi=\int_{[0,1]^N}|\nabla\chi|^2>0$), while $\int(w\chi(k\cdot))^q\,dx \to \bar{d}_\chi\int w^q\,dx>0$ (where $\bar{d}_\chi=\int_{[0,1]^N}\chi^q>0$), so $Q(1,k) > \beta$ for $k$ large. The intermediate value theorem gives $t_k$ with $Q(t_k,k)=\beta$; thus, we can take $v = v_{t_k,k}$.
\end{proof}

\begin{prop}\label{CurveIsContinuousAndDecreasing}
    The mapping $A : \beta \mapsto\alpha_\beta$ is continuous and strictly decreasing on $(\lambda_1(\Omega),+\infty)$. Moreover, $A(A(\beta)) = \beta$ for all $\beta>\lambda_1(\Omega)$. In particular, the first non-trivial curve is symmetric with respect to the diagonal.
\end{prop}

\begin{proof}

We first show that the map $\beta \mapsto\alpha_\beta$ is strictly decreasing on $(\lambda_1(\Omega),+\infty)$. Let $\beta_1,\beta_2 > \lambda_1(\Omega)$ such that $\beta_1<\beta_2$. Since $\int_\Omega |\nabla \bar{u}_{\beta_1}^-|^2\,dx = \alpha_{\beta_1}$ and $\int_\Omega |\nabla \bar{u}_{\beta_1}^+|^2\,dx = \beta_1 < \beta_2$, then by Lemma \ref{ModificationOfSupports}, we can construct a function $v \in H_0^1(\Omega)$ such that $\norm{v^+}_{L^q(\Omega)} = \norm{v^-}_{L^q(\Omega)} = 1$, $\int_\Omega |\nabla v^+|^2\,dx < \beta_2$, and $\int_\Omega |\nabla v^-|^2\,dx < \alpha_{\beta_1}$. Moreover, since $\int_\Omega |\nabla v^+|^2\,dx < \beta_2$ by Lemma \ref{DecreasingTheNorm}, we can construct a function $w \in H_0^1(\Omega)$ such that $0 \leq w \leq v^+$ and $\int_\Omega |\nabla w|^2\,dx = \beta_2 \norm{w}_{L^q(\Omega)}^2$. We then set $\hat{u}_{\beta_2} = \frac{w}{\norm{w}_{L^q(\Omega)}} - v^-$; we have $\norm{\hat{u}_{\beta_2}^+}_{L^q(\Omega)} = \norm{\hat{u}_{\beta_2}^-}_{L^q(\Omega)} = 1$, $\int_\Omega |\nabla \hat{u}_{\beta_2}^+|^2\,dx = \beta_2$, and $\int_\Omega |\nabla \hat{u}_{\beta_2}^-|^2\,dx = \int_\Omega |\nabla v^-|^2\,dx < \alpha_{\beta_1}$. Therefore, the function $\hat{u}_{\beta_2}$ is admissible
$$
    \alpha_{\beta_2} \leq \int_\Omega |\nabla \hat{u}_{\beta_2}^-|^2\,dx < \alpha_{\beta_1}.
$$
Therefore, the mapping is strictly decreasing, as announced.

We now show that $A(A(\beta))=\beta$ for all $\beta>\lambda_1(\Omega)$. Since $(A(\beta),\beta)\in\Sigma$ and $\Sigma$ is symmetric with respect to the diagonal, $(\beta,A(\beta))\in\Sigma$; Proposition \ref{SmallestFucikPair}, applied with $A(\beta)$ in place of $\beta$, gives $A(A(\beta))\leq\beta$. If the inequality were strict, then $A(A(A(\beta)))>A(\beta)$, since $A$ is strictly decreasing. On the other hand, $(A(A(\beta)),A(\beta))\in\Sigma$, hence $(A(\beta),A(A(\beta)))\in\Sigma$, and Proposition \ref{SmallestFucikPair}, applied with $A(A(\beta))$ in place of $\beta$, gives $A(A(A(\beta)))\leq A(\beta)$, a contradiction. Hence $A(A(\beta))=\beta$. In particular, $A$ maps $(\lambda_1(\Omega),+\infty)$ onto itself; being strictly decreasing, it has no jump discontinuities, and is therefore continuous.
\end{proof}

\begin{prop}\label{CurveAsymptotes}
    For all $\beta>\lambda_1(\Omega)$, let $\bar{u}_\beta$ be the minimiser in \eqref{Definition_alpha_beta} given by Theorem \ref{FucikSelectedEigenfunction}, so that $\norm{\bar u_\beta^+}_{L^q(\Omega)}=\norm{\bar u_\beta^-}_{L^q(\Omega)}=1$, $\int_\Omega|\nabla\bar u_\beta^+|^2\,dx=\beta$ and $\int_\Omega|\nabla\bar u_\beta^-|^2\,dx=\alpha_\beta$. Let $e_1$ be the first (positive) $(2,q)$-eigenfunction of the Dirichlet Laplacian, with $\norm{e_1}_{L^q(\Omega)}=1$.
When $\beta\rightarrow+\infty$, we have $\alpha_\beta \rightarrow \lambda_1(\Omega)$, and $\bar{u}_\beta^- \rightarrow e_1$ in $H_0^1(\Omega)$. When $\beta\rightarrow\lambda_1(\Omega)$, we have $\alpha_\beta \rightarrow +\infty$, and $\bar{u}_\beta^+ \rightarrow e_1$ in $H_0^1(\Omega)$.
\end{prop}

\begin{proof}
Since the application $\beta\mapsto\alpha_\beta$ is decreasing, we have for all $\beta>\lambda_1(\Omega)$,
\[
\begin{split}
    \lambda_1(\Omega) < \alpha_\beta &= \inf\left\{\int_\Omega |\nabla u^-|^2\,dx \,;\, u\in H_0^1(\Omega), \norm{u^\pm}_{L^q(\Omega)}=1 , \int_\Omega |\nabla u^+|^2\,dx=\beta \right\} \\
    & = \inf\left\{\int_\Omega |\nabla u^-|^2\,dx \,;\, u\in H_0^1(\Omega), \norm{u^\pm}_{L^q(\Omega)}=1 , \int_\Omega |\nabla u^+|^2\,dx\leq\beta \right\}.
\end{split}
\]
By letting $\beta$ go to $+\infty$, we obtain
\[
\begin{split}
    \lambda_1(\Omega) &\leq \lim_{\beta\rightarrow+\infty} \alpha_\beta = \inf\{\alpha_\beta \,;\, \beta>\lambda_1(\Omega)\} \\
    & \leq \inf\left\{\int_\Omega |\nabla u^-|^2\,dx \,;\, u\in H_0^1(\Omega), \norm{u^+}_{L^q(\Omega)} = \norm{u^-}_{L^q(\Omega)} = 1 \right\} \\
    & \leq \int_{\Omega_r} |\nabla e_1(\Omega_r)|^2\,dx \,\,\,\, , \,\,\,\, \forall r\in \left(0, \frac{1}{2} diam(\Omega)\right),
\end{split}
\]
where for all $r\in (0,\frac{1}{2} diam(\Omega))$, $\Omega_r=\Omega\setminus \bar{B}(x_0,r)$ (with $x_0\in\partial\Omega$) and $e_1(\Omega_r)$ is the eigenfuntion associated to the $(2,q)$-eigenvalue $\lambda_1(\Omega_r)$, normalised in $L^q$.
By letting $r$ go to $0$, we obtain that $\int_{\Omega_r} |\nabla e_1(\Omega_r)|^2\,dx = \lambda_1(\Omega_r) \rightarrow \lambda_1(\Omega)$, by Lemma \ref{DecreasingSequenceOfLambda1}. Therefore, $\lim_{\beta\rightarrow+\infty} \alpha_\beta = \lambda_1(\Omega)$.

For all $\beta>\lambda_1(\Omega)$, $\norm{\bar{u}_\beta^-}_{L^q(\Omega)} = 1$ and $\int_\Omega |\nabla \bar{u}_\beta^-|^2\,dx = \alpha_\beta$. Since $\lim_{\beta\rightarrow+\infty} \alpha_\beta = \lambda_1(\Omega)$, the sequence $(\bar{u}_\beta^-)_\beta$ converges to a function $\bar{u}\in H_0^1(\Omega)$ strongly in $L^q(\Omega)$, weakly in $H_0^1(\Omega)$, and a.e. in $\Omega$. Hence, by the definition of $\lambda_1(\Omega)$, we obtain
$$
    \lambda_1(\Omega) =\lambda_1(\Omega) \left[\int_\Omega |\bar{u}|^q\,dx\right]^{\frac{2}{q}} \leq \int_\Omega |\nabla \bar{u}|^2\,dx \leq \lim_{\beta\rightarrow+\infty} \int_\Omega |\nabla \bar{u}_\beta^-|^2\,dx = \lim_{\beta\rightarrow+\infty} \alpha_\beta = \lambda_1(\Omega).
$$
Therefore, $(\bar{u}_\beta^-)_\beta$ converges strongly in $H_0^1(\Omega)$ to $\bar{u}$ (up to a subsequence, then the full family by uniqueness of the limit). Moreover, $\bar{u} \geq 0$, $\norm{\bar{u}}_{L^q(\Omega)} = 1$, and $\int_\Omega |\nabla \bar{u}|^2 \,dx = \lambda_1(\Omega)$; hence, $\bar{u}$ corresponds to the infimum in \eqref{FirstEV}. By Propositions \ref{FirstEigenFunctionsArePositive} and \ref{PositiveEigenfunctionsAreFirst}, $\bar{u}$ corresponds to the first $(2,q)$-eigenfunction $e_1$.

We now study the limit $\beta\to\lambda_1(\Omega)^+$. We assume by contradiction that there exists a decreasing sequence $(\beta_n)_n$ that converges to $\lambda_1(\Omega)$, such that $\sup_n\alpha_{\beta_n}<+\infty$.
For all $n\in\N$, the minimiser $\bar{u}_{\beta_n}$ given by Proposition \ref{PDE_limit} satisfies $\norm{\bar{u}_{\beta_n}^+}_{L^q(\Omega)} = \norm{\bar{u}_{\beta_n}^-}_{L^q(\Omega)} = 1$, $\int_\Omega|\nabla \bar{u}_{\beta_n}^+|^2 \,dx = \beta_n$, and $\int_\Omega |\nabla \bar{u}_{\beta_n}^-|^2 \,dx = \alpha_{\beta_n}$; hence, the sequence $(\bar{u}_{\beta_n})_n$ is bounded in $H_0^1(\Omega)$. Up to a subsequence, the sequence $\bar{u}_{\beta_n}$ converges to $w$ strongly in $L^q(\Omega)$, weakly in $H_0^1(\Omega)$, and a.e. in $\Omega$. Thus, $\norm{w^+}_{L^q(\Omega)} = \norm{w^-}_{L^q(\Omega)} = 1$ and,
by weak lower semi-continuity, $\int_\Omega |\nabla w^+|^2\,dx \leq \liminf_n \beta_n = \lambda_1(\Omega)$. Hence, $w^+$ attains the infimum in \eqref{FirstEV}; by Propositions \ref{FirstEigenFunctionsArePositive} and~\ref{PositiveEigenfunctionsAreFirst}, $w^+ = e_1 > 0$ a.e. in $\Omega$, whence $w^- = 0$ a.e., contradicting $\norm{w^-}_{L^q(\Omega)} = 1$. Therefore, $\lim_{\beta\to\lambda_1(\Omega)^+} \alpha_\beta = +\infty$.

The same argument, together with the convergence of the norms
$\displaystyle
    \int_\Omega |\nabla \bar{u}_\beta^+|^2 \,dx = \beta \to \lambda_1(\Omega) = \int_\Omega |\nabla e_1|^2 \,dx,
$
shows that $\bar{u}_\beta^+$ converges to $e_1$ in $H_0^1(\Omega)$, as $\beta$ goes to $\lambda_1(\Omega)$.

\end{proof}

All those properties of the \textbf{first curve} in the Fu\v{c}\'\i k $(2,q)$-spectrum show us that it has a unique diagonal point. The continuous curve $\{(\alpha,\beta)\in\Sigma ; \beta>\lambda_1(\Omega),\, \alpha = \alpha_\beta\}$ intersects the line $\{(\alpha,\beta)\in\R^2 ; \alpha = \beta>\lambda_1(\Omega)\}$ in a unique point.
Indeed, the mapping $g\,:\, \beta \mapsto \alpha_\beta - \beta$ is continuous and strictly decreasing on $(\lambda_1(\Omega),+\infty)$ by Proposition \ref{CurveIsContinuousAndDecreasing}, with $\lim_{\beta\rightarrow \lambda_1(\Omega)^+} g(\beta) = +\infty$ and $\lim_{\beta\rightarrow +\infty} g(\beta) = -\infty$ by Proposition \ref{CurveAsymptotes}; hence, there exists a unique $\gamma > \lambda_1(\Omega)$ with $\alpha_\gamma = \gamma$.
\begin{defn}\label{defgamma}
{We denote by $\gamma=\gamma(\Omega)>\lambda_1(\Omega)$ the unique value such that $\alpha_{\gamma(\Omega)}=\gamma(\Omega)$; the point $(\gamma(\Omega),\gamma(\Omega))$ is the unique intersection of the first non-trivial curve with the diagonal.}
{We also define the Fu\v{c}\'\i k-selected second value by $\lambda_2^F(\Omega) := 2^{1 - \frac{2}{q}}\gamma(\Omega)$, independently of the nodal balance of the selected solution.}
Let $u_\gamma \in H_0^1(\Omega)$ be the eigenfunction associated to the pair $(\gamma,\gamma)$, as given by Proposition \ref{PDE_limit}; thus, by setting $t_\gamma = \norm{u_\gamma^+}_{L^q(\Omega)} > 0$, $u_\gamma$ solves
\begin{equation}\label{FucikDiagonalEquation}
    -\Delta u_\gamma = -\gamma (u_\gamma^-)^{q-1} + \gamma t_\gamma^{2-q} (u_\gamma^+)^{q-1}.
\end{equation}
\end{defn}
In the following, we shall call $u_\gamma$ a \textbf{diagonal Fu\v{c}\'\i k-selected solution}. We notice that we have no information concerning the value of $t_\gamma$.

\begin{prop}\label{SecondFucikEigenvalue}
    The function $u_\gamma$ is a classical $(2,q)$-eigenfunction, namely there exists a $(2,q)$-eigenvalue $\lambda\in\sigma$ such that
    $$
        -\Delta u_\gamma = \lambda \norm{u_\gamma}_{L^q(\Omega)}^{2-q} |u_\gamma|^{q-2} u_\gamma,
    $$
    if and only if $t_\gamma = 1$. In that case, $\lambda$ equals the Fu\v{c}\'\i k-selected second value defined above, and we have
    $
      \lambda=  \lambda_2^F(\Omega) = 2^{1 - \frac{2}{q}} \gamma.
    $
\end{prop}

\begin{proof}
Compare \eqref{FucikDiagonalEquation} with \eqref{Eigen2q}: if $t_\gamma=1$ then $\norm{u_\gamma}_{L^q(\Omega)}^q=2$ and $-\Delta u_\gamma=\gamma|u_\gamma|^{q-2}u_\gamma=2^{1-\frac2q}\gamma\norm{u_\gamma}_{L^q(\Omega)}^{2-q}|u_\gamma|^{q-2}u_\gamma$. Conversely, the two equations coincide only if $t_\gamma^{2-q}=1$.
\end{proof}

\begin{rem}
    Hence, the first curve in the Fu\v{c}\'\i k $(2,q)$-spectrum and its unique diagonal point $(\gamma,\gamma)$ are fully defined, but the identification of $u_\gamma$ with a classical $(2,q)$-eigenfunction at the value $\lambda_2^F(\Omega) = 2^{1 - \frac{2}{q}} \gamma$ requires the additional balance condition $\norm{u_\gamma^+}_{L^q(\Omega)} = \norm{u_\gamma^-}_{L^q(\Omega)}$. This condition disappears when $q=2$, since the norms in the Fu\v{c}\'\i k equation vanish, and the two coefficients coincide independently of the relative norms of $u_\gamma^+$ and $u_\gamma^-$ (as explained in Remark \ref{Remark q=2}). This shows why the sublinear case $1<q<2$ is more delicate than the linear case $q=2$. We prove in Proposition \ref{lambda1_less_lambda2F} below that, in any case, $\lambda_2^F(\Omega)>\lambda_1(\Omega)$. The balance condition $t_\gamma=1$ is studied in \cite{PartII}, in the limit $q\to2^-$.

\end{rem}

\subsection{Regularity of the Fu\v{c}\'\i k $(2,q)$-eigenfunctions}

Before concluding this section, we state a theorem regarding the regularity of \textit{all} solutions of the Fu\v{c}\'\i k equation \eqref{Fucik2q}. We may notice that the result is similar to the regularity result of the classical $(2,q)$-eigenfunctions, given in Propositions \ref{InteriorRegularity} and \ref{BoundaryRegularity}. We shall use the following lemma, which is a particular case of \cite[Theorems 4.1 and 4.2]{Stampacchia} (where the right-hand side is written in divergence form, $b=\operatorname{div}(\nabla w)$ with $w$ the Newtonian potential of $b$; see also \cite[Theorem 8.15]{Gilbarg_Trudinger}).

\begin{prop}\label{StampacchiaLemma}
    Let $u\in H_0^1(\Omega)$ be a solution of the equation
    $$
        -\Delta u = au + b,
    $$
    where $a,b \in L^\infty(\Omega)$. Therefore, there exist constants $K = K(N,\norm{a}_{L^\infty(\Omega)}) > 0$ and $C = C(N,\Omega,\norm{a}_{L^\infty(\Omega)}) > 0$ such that
    $$
        \norm{u}_{L^\infty(\Omega)} \leq K \norm{b}_{L^\infty(\Omega)} |\Omega|^{\frac{2}{N}} + C\norm{u}_{L^2(\Omega)}.
    $$
    In particular, we obtain that $u\in L^\infty(\Omega)$.
\end{prop}

\begin{prop}\label{FucikRegularity}
    Let $(\alpha,\beta)\in\Sigma$, and let $u$ be an associated Fu\v{c}\'\i k $(2,q)$-eigenfunction. Then, $u\in L^\infty(\Omega) \cap C_{loc}^{2,q-1}(\Omega)$. If we assume that $\partial\Omega$ is of class $C^{2,q-1}$, then $u\in C^{2,q-1}(\overline\Omega)$.
\end{prop}

\begin{proof}
The eigenfunction $u\in H_0^1(\Omega)$ satisfies the following problem:
$$
    -\Delta u = f, \qquad \textnormal{ with } f = -\alpha\norm{u^-}_{L^q(\Omega)}^{2-q}(u^-)^{q-1} + \beta\norm{u^+}_{L^q(\Omega)}^{2-q}(u^+)^{q-1}.
$$
We first notice that there exists a constant $C>0$ such that $|f| \leq C |u|^{q-1}$.

We first prove that $u\in L^\infty(\Omega)$. We define $a = \frac{f}{u} \chi_{\{ |u| > 1 \}}$ and $b = f \chi_{\{ |u| \leq 1 \}}$, such that $f = au+b$. Since $|f| \leq C |u|^{q-1}$ and $q<2$, we obtain that $|a| \leq C$. Moreover, we have $|b|\leq C$. Therefore, $a,b \in L^\infty(\Omega)$. The equation satisfied by $u\in H_0^1(\Omega)$ can thus be written as
$$
    -\Delta u - au = b.
$$
Hence, applying Proposition \ref{StampacchiaLemma} yields $u\in L^\infty(\Omega)$.

Since $\Omega$ is bounded, we obtain that $f\in L^r(\Omega)$ for all $r\in (1,+\infty)$. The interior Calder\'on--Zygmund estimate from Proposition \ref{RegularityTheorems} gives that $u\in W^{2,r}_{loc}(\Omega)$, for all $r\in (1,+\infty)$. By Morrey's embedding, we obtain that, by choosing $r>N$, $u\in C^{1,\alpha}_{loc}(\Omega)$ with $\alpha=1-\frac{N}{r}$, and thus $u\in C^{0,1}_{loc}(\Omega)$. Since the functions $s\mapsto(s^+)^{q-1}$ and $s\mapsto(s^-)^{q-1}$ are $(q-1)$-H\"older continuous on $\R$, we obtain $f\in C^{0,q-1}_{loc}(\Omega)$. The interior Schauder estimate from Proposition \ref{RegularityTheorems} implies that $u \in C^{2,q-1}_{loc}(\Omega)$.

To obtain the global regularity result if $\partial\Omega$ is of class $C^{2,q-1}$, it suffices to repeat the previous argument and use the global regularity theorems given in Proposition \ref{BoundaryRegularityTheorems}. Since $\Omega$ is bounded, we obtain that $f\in L^r(\Omega)$ for all $r\in (1,+\infty)$. The global Calder\'on--Zygmund estimate gives that $u\in W^{2,r}(\Omega)$, for all $r\in (1,+\infty)$. By Morrey's embedding, we obtain that, by choosing $r>N$, $u\in C^{1,\alpha}(\overline{\Omega})$ with $\alpha=1-\frac{N}{r}$. Therefore, $u\in C^{0,1}(\overline{\Omega})$, and thus $f\in C^{0,q-1}(\overline{\Omega})$. The global Schauder estimate then implies that $u \in C^{2,q-1}(\overline{\Omega})$.
\end{proof}

\begin{prop}\label{lambda1_less_lambda2F}
    Let $\Omega$ be a bounded connected open set. Then $\lambda_2^F(\Omega)=2^{1-\frac2q}\gamma(\Omega)>\lambda_1(\Omega)$.
\end{prop}

\begin{proof}
Let $\bar u_\gamma=-u_\gamma^-+t_\gamma^{-1}u_\gamma^+$, so that $\norm{\bar u_\gamma^+}_{L^q(\Omega)}=\norm{\bar u_\gamma^-}_{L^q(\Omega)}=1$ and, by Theorem \ref{FucikSelectedEigenfunction} with $\alpha_\gamma=\gamma$, $\int_\Omega|\nabla\bar u_\gamma^+|^2\,dx=\int_\Omega|\nabla\bar u_\gamma^-|^2\,dx=\gamma$. The function $w=|\bar u_\gamma|=\bar u_\gamma^++\bar u_\gamma^-$ satisfies $\int_\Omega|\nabla w|^2\,dx=2\gamma$ and $\norm{w}_{L^q(\Omega)}^q=2$, so that \eqref{FirstEV} gives $\lambda_1(\Omega)\leq 2\gamma\,2^{-\frac2q}=\lambda_2^F(\Omega)$. Suppose that equality holds. Then $w$ is a first $(2,q)$-eigenfunction: it is continuous in $\Omega$ by Proposition \ref{InteriorRegularity}, non-negative and superharmonic, hence $w>0$ in $\Omega$ by the strong maximum principle.
On the other hand, $u_\gamma$ is continuous in $\Omega$ by Proposition \ref{FucikRegularity}, hence so is $\bar u_\gamma=-u_\gamma^-+t_\gamma^{-1}u_\gamma^+$; since $u_\gamma$ changes sign (Proposition \ref{PDE_limit}), so does $\bar u_\gamma$, and, $\Omega$ being connected, it vanishes at some point of $\Omega$; this contradicts $w>0$.
\end{proof}

\section{A Hong--Krahn--Szeg\H{o} inequality for the Fu\v{c}\'\i k-selected second value}

We have defined in Definition \ref{defgamma} the value $\gamma > \lambda_1(\Omega)$ as the unique value such that $\alpha_{\gamma} = \gamma$; it corresponds to the diagonal point of the first non-trivial Fu\v{c}\'\i k curve. We shall denote in this section $\gamma = \gamma(\Omega)$ to show the dependence of this value on the connected bounded open set $\Omega$. We can define the value
\begin{equation}\label{lambda_2^F}
    \lambda_2^F(\Omega) := 2^{1 - \frac{2}{q}}\,\gamma(\Omega).
\end{equation}
As seen in Proposition \ref{SecondFucikEigenvalue}, $\lambda_2^F(\Omega)$ is a $(2,q)$-eigenvalue, with eigenfunction $u_\gamma$, when the balance condition $t_{\gamma(\Omega)}=1$ holds, and $\lambda_2^F(\Omega)>\lambda_1(\Omega)$ by Proposition \ref{lambda1_less_lambda2F}. In any case, it admits a sharp shape optimisation principle.

Since $\gamma(\Omega) = \alpha_{\gamma(\Omega)}$, Theorem \ref{FucikSelectedEigenfunction} gives the following variational characterisation
\begin{equation}\label{Definition_gamma}
    \gamma(\Omega) = \inf\left\{\int_\Omega |\nabla u^-|^2\,dx \,;\, u\in H_0^1(\Omega), \norm{u^+}_{L^q(\Omega)} = \norm{u^-}_{L^q(\Omega)} = 1 , \int_\Omega |\nabla u^+|^2\,dx=\gamma(\Omega) \right\}.
\end{equation}
Since this definition is self-referencing, we would like to obtain a simpler variational characterisation of $\gamma(\Omega)$. We shall relate it to a two-phase characterisation as follows.

We set
$$
    E(u) = \int_\Omega |\nabla u|^2\,dx,
$$
and we define the two-phase class
$$
    \mathcal{A}_2(\Omega) = \left\{(v,w)\in H_0^1(\Omega)^2\,;\, v,w \geq 0, \, vw=0\ \text{a.e. in }\Omega, \, \norm{v}_{L^q(\Omega)} = \norm{w}_{L^q(\Omega)} = 1 \right\}.
$$
We define
\begin{equation}\label{Definition_Gamma}
    \Gamma(\Omega) = \inf_{(v,w)\in\mathcal A_2(\Omega)} \max \{E(v),E(w)\}.
\end{equation}

\begin{prop}\label{Gamma = gamma}
    Let $\Omega$ be a bounded open connected set of $\R^N$. Then
    $
        \Gamma(\Omega)=\gamma(\Omega).
    $
    Moreover, the infimum in \eqref{Definition_Gamma} is attained.
\end{prop}

\begin{proof}
Let $u_\gamma$ be a minimiser in the definition \eqref{Definition_alpha_beta} of $\alpha_\gamma$, such that $\norm{u_\gamma^+}_{L^q(\Omega)} = \norm{u_\gamma^-}_{L^q(\Omega)} = 1$. We have $(u_\gamma^+, u_\gamma^-) \in \mathcal{A}_2(\Omega)$. Since $\alpha_\gamma = \gamma$, we have
$
    E(u_\gamma^+) = E(u_\gamma^-) = \gamma(\Omega).
$
Therefore, $\Gamma(\Omega) \leq \gamma(\Omega)$.

Conversely, let $(v,w)\in\mathcal{A}_2(\Omega)$. We set $\tau = \max\{E(v),E(w)\}$, and thus $\Gamma(\Omega) = \inf_{(v,w)\in\mathcal A_2(\Omega)} \tau$. Since $\Omega$ is connected, and the two phases are non-trivial ($\norm{v}_{L^q(\Omega)} = \norm{w}_{L^q(\Omega)} = 1$) and disjoint ($vw = 0$), neither phase can be a first $(2,q)$-eigenfunction on $\Omega$.
Hence, $E(v),E(w)>\lambda_1(\Omega)$, and thus, $\tau>\lambda_1(\Omega)$.
Assume for example, that $E(v)\leq E(w)=\tau$. By Lemma \ref{DecreasingTheNorm}, $v$ can be replaced by a function $\bar{v}$ such that $0 \leq \bar{v} \leq v$ and $E(\bar{v}) = \tau \norm{\bar{v}}^{2}_{L^q(\Omega)}$.

Setting $\tilde{v} = \frac{\bar{v}}{\norm{\bar{v}}_{L^q(\Omega)}}$ produces a function supported in $Supp(v)$ such that $\norm{\tilde{v}}_{L^q(\Omega)} = 1$ and $E(\tilde{v}) = \tau$. Hence, there exist non-negative functions $\tilde v,\tilde w\in H_0^1(\Omega)$, with disjoint supports ($vw = 0$) and
unit $L^q$-norms, such that
$
    E(\tilde v) = E(\tilde w) = \tau.
$
Using $\tilde v - \tilde w$ as an admissible function in the definition of $\alpha_\tau$ gives $\alpha_\tau\leq \tau$. The map $\beta \mapsto \alpha_\beta-\beta$ is strictly decreasing and admits $\gamma(\Omega)$ as its unique zero. Therefore, since $\alpha_\tau - \tau \leq 0$, we have $\tau \geq \gamma(\Omega)$. Taking the infimum over $\mathcal A_2(\Omega)$ yields $\Gamma(\Omega) \geq \gamma(\Omega)$.

\end{proof}

\begin{rem}\label{rem:partition}
{The two-phase characterisation of Proposition \ref{Gamma = gamma} can be rephrased as an optimal partition problem. We refer to \cite{Bucur_Buttazzo} for the notions of capacity, quasi-open sets and quasi-continuous representatives. For a quasi-open set $\omega\subset\Omega$, let $H_0^1(\omega)$ be the closed subspace of the functions $u\in H_0^1(\Omega)$ whose quasi-continuous representative vanishes quasi-everywhere in $\Omega\setminus\omega$, and let $\lambda_1(\omega)$ be defined by \eqref{FirstEV} with $H_0^1(\omega)$ in place of $H_0^1(\Omega)$ (with $\lambda_1(\omega)=+\infty$ if $H_0^1(\omega)=\{0\}$); by the compactness of the embedding $H_0^1(\Omega)\hookrightarrow L^q(\Omega)$, the infimum is attained, and it is attained by a non-negative function, since $|u|\in H_0^1(\omega)$ whenever $u\in H_0^1(\omega)$. Then
$$
\Gamma(\Omega)=\inf\Big\{\max\{\lambda_1(\omega_1),\lambda_1(\omega_2)\}\,;\ \omega_1,\omega_2\subset\Omega\ \text{quasi-open},\ |\omega_1\cap\omega_2|=0\Big\}.
$$
Indeed, from one side, if $\omega_1,\omega_2$ are quasi-open with $|\omega_1\cap\omega_2|=0$ and $\lambda_1(\omega_1),\lambda_1(\omega_2)<+\infty$, the non-negative $L^q$-normalised first eigenfunctions $v\in H_0^1(\omega_1)$, $w\in H_0^1(\omega_2)$ satisfy $vw=0$ a.e. in $\Omega$ (outside the null set $\omega_1\cap\omega_2$, at least one of them vanishes quasi-everywhere), so that $(v,w)\in\mathcal A_2(\Omega)$ with $E(v)=\lambda_1(\omega_1)$ and $E(w)=\lambda_1(\omega_2)$.

On the other side, if $(v,w)\in\mathcal A_2(\Omega)$, the sets $\omega_1=\{v>0\}$ and $\omega_2=\{w>0\}$, defined through quasi-continuous representatives, are quasi-open, $|\omega_1\cap\omega_2|=0$ since $vw=0$ a.e., and $v\in H_0^1(\omega_1)$, $w\in H_0^1(\omega_2)$, whence $\lambda_1(\omega_1)\leq E(v)$ and $\lambda_1(\omega_2)\leq E(w)$. In the linear case $q=2$, the characterisation of the first non-trivial Fu\v{c}\'\i k curve through optimal partition problems is due to Conti, Terracini and Verzini \cite[Theorem 1.2]{Conti_Terracini_Verzini}.}
\end{rem}

For $m>0$, we denote by $B_m$ any ball of measure $m$. {The following inequality is classical; we include the proof for the sake of completeness, since the equality case is used in Theorem \ref{HKS_lambda2F}.}

\begin{Lemma}\label{RearrangingThePositivePart}
    Let $z\in H_0^1(\Omega)$ be a non-negative and non-trivial function, and set $m_z = |\{z>0\}|$. Then
    \begin{equation}\label{RearrangingThePositivePart_eq}
        \frac{E(z)}{\norm{z}_{L^q(\Omega)}^2} = \frac{\int_\Omega |\nabla z|^2 \,dx}{\norm{z}_{L^q(\Omega)}^2} \geq \lambda_1(B_{m_z}).
    \end{equation}
    If equality holds, then, up to a translation and modification on a null set, $\{z>0\}$ is a ball of measure $m_z$, and $z$ is a positive first $(2,q)$-eigenfunction on that ball.
\end{Lemma}

\begin{proof}
Let $z^*$ be the Schwarz rearrangement of $z$. Then, $z^*\in H_0^1(B_{m_z})$, Cavalieri's principle gives $\norm{z}_{L^q(\Omega)} = \norm{z^*}_{L^q(B_{m_z})}$, and the P\'olya--Szeg\H{o} inequality gives $E(z^*)\leq E(z)$. Hence
$$
    \frac{E(z)}{\norm{z}_{L^q(\Omega)}^2} \geq \frac{E(z^*)}{\norm{z^*}_{L^q(B_{m_z})}^2} \geq \lambda_1(B_{m_z}),
$$
the last inequality following from the definition \eqref{FirstEV} of $\lambda_1(B_{m_z})$. This yields the desired inequality \eqref{RearrangingThePositivePart_eq}.

If equality in \eqref{RearrangingThePositivePart_eq} holds, then
$$
    \frac{E(z)}{\norm{z}_{L^q(\Omega)}^2} = \frac{E(z^*)}{\norm{z^*}_{L^q(B_{m_z})}^2} = \lambda_1(B_{m_z}),
$$
Thus, $z^*$ is a first eigenfunction on $B_{m_z}$, and there is equality in the P\'olya--Szeg\H{o} inequality. Every positive first $(2,q)$-eigenfunction $u$ on a ball $B$ of radius $R$ is radial. Moreover, $u$ must be strictly decreasing in the radial variable, because the radial profile $U$ defined by $U(r) = u(x)$, where $r = |x|$, satisfies the following ordinary differential equation
$$
    -\left(r^{N-1} U'(r)\right)' = \lambda_1(B) \norm{u}_{L^q(B)}^{2-q} r^{N-1} U(r)^{q-1}, \qquad U'(0) = 0, \quad U(R) = 0.
$$
Integrating from $0$ to $r$, we obtain
$$
    r^{N-1} U'(r) = - \lambda_1(B) \norm{u}_{L^q(B)}^{2-q} \int_0^r s^{N-1} U(s)^{q-1} \,ds < 0,
$$
for $0<r<R$. Thus, $U$ is strictly decreasing. In particular, the set on which its radial derivative vanishes has measure zero away from the origin. Applying this to $u=z^*$, we obtain
$$
    \left| \{\nabla z^* = 0\} \cap \{0 < z^* < \operatorname*{ess\,sup} z^*\}\right| = 0.
$$
\cite[Theorem 1.1]{Brothers_Ziemer} therefore applies and yields that $z$ is a translate of $z^*$, up to a set of zero measure. Hence, $\{z>0\}$ is a ball of measure $m_z$ up to a set of measure zero, and $z$ is its positive first $(2,q)$-eigenfunction.
\end{proof}

We can now prove a shape optimisation theorem of $\gamma(\Omega)$ among all connected bounded open subsets $\Omega$ of $\R^N$. We will prove that the infimum of $\gamma(\Omega)$ is not attained among the considered class of subsets. This result is similar to that one obtained for the second eigenvalue of the fractional $p$-Laplacian as done in \cite{Brasco_Parini}.

\begin{thm}\label{HKS_lambda2F}
    Let $\Omega$ be a bounded open connected set of $\R^N$, and set $V=|\Omega|$. Then
    \begin{equation}\label{HKS_gamma}
        \gamma(\Omega) \geq \lambda_1(B_{V/2}),
    \end{equation}
    and consequently
    \begin{equation}\label{HKS_lambda2F_eq}
        \lambda_{2}^F(\Omega) = 2^{1 - \frac{2}{q}}\gamma(\Omega) \geq 2^{1 - \frac{2}{q}}\lambda_1(B_{V/2}).
    \end{equation}
    For $N\geq 2$, both inequalities are strict in the class of connected open sets. Moreover, they are sharp in the following sense:
    \begin{equation}\label{ConnectedShapeInfimum}
        \inf \left\{\lambda_{2}^F(\Omega); \Omega\subset\R^N\text{ bounded, open and connected}, |\Omega| = V \right\} = 2^{1 - \frac{2}{q}}\lambda_1(B_{V/2}).
    \end{equation}
    For $N \geq 2$, the infimum is not attained; for $N=1$, it is attained by every interval of length $V$.
\end{thm}

\begin{proof}
Let $(v,w)\in\mathcal{A}_2(\Omega)$, and set $m_v = |\{v>0\}|$ and $m_w = |\{w>0\}|$. Since $vw=0$ a.e. in $\Omega$, we have $m_v + m_w \leq |\Omega|=V$. Lemma \ref{RearrangingThePositivePart}, together with $\norm{v}_{L^q(\Omega)} = \norm{w}_{L^q(\Omega)} = 1$, yields
$$
    E(v)\geq\lambda_1(B_{m_v}), \qquad E(w)\geq\lambda_1(B_{m_w}).
$$
Since $m_v + m_w \leq |\Omega|=V$, at least one of $m_v$ or $m_w$ is not larger than $V/2$.  Since $m\to \lambda_1(B_m)$ is strictly decreasing by Corollary \ref{StrictMonotonicityEV1}, we obtain
$
    \max\{E(v),E(w)\}\geq\lambda_1(B_{V/2}).
$
Taking the infimum over $\mathcal{A}_2(\Omega)$ and using Proposition \ref{Gamma = gamma} prove the inequality
$
    \gamma(\Omega) \geq \lambda_1(B_{V/2}).
$
From the definition \eqref{lambda_2^F} of $\lambda_2^F(\Omega)$, it  follows that
$
    \lambda_{2}^F(\Omega) \geq 2^{1 - \frac{2}{q}}\lambda_1(B_{V/2}).
$

Next, we exclude equality when $\Omega$ is connected for $N\geq 2$. Suppose that $\gamma(\Omega) = \lambda_1(B_{V/2})$. By Proposition \ref{Gamma = gamma}, there exists a pair $(v,w)\in\mathcal{A}_2(\Omega)$ such that $\gamma(\Omega) = \max\{E(v),E(w)\}$. The previous inequalities then imply
$$
    m_v = m_w = \frac{V}{2}, \qquad E(v) = E(w) = \lambda_1(B_{V/2}).
$$
Indeed, if one of the two measures were strictly smaller than $\frac{V}{2}$, the strict monotonicity of $m\mapsto\lambda_1(B_m)$ would give a strict lower bound for the corresponding energy, and $\gamma(\Omega) < \max\{E(v),E(w)\}$, which yields a contradiction. Lemma \ref{RearrangingThePositivePart} therefore shows that, up to zero-measure sets, $\{v>0\}$ and $\{w>0\}$ are two balls, denoted respectively by $B_1$ and $B_2$, of equal measure $\frac{V}{2}$. Since $vw=0$ a.e. in $\Omega$, $B_1$ and $B_2$ are disjoint. Moreover,
$$
    |\Omega\setminus(B_1\cup B_2)| = V-|B_1|-|B_2| = 0.
$$
{We claim that this is impossible for a connected open set $\Omega$ for $N\geq 2$.} {Indeed, $\Omega\subset B_1\cup B_2$. Otherwise, a point $x\in\Omega\setminus(B_1\cup B_2)$ has a neighbourhood contained in $\Omega$ whose intersection with the complement of $B_1\cup B_2$ has positive measure: this is immediate if $x\notin\overline{B_1\cup B_2}$, and also holds at a boundary point of two disjoint balls when $N\geq2$. This contradicts $|\Omega\setminus(B_1\cup B_2)|=0$. Thus, up to a set of zero measure}
\[
{\Omega= B_1\cup B_2.}
\]
We have thus obtained a decomposition of $\Omega$ into two disjoint, non-empty, open subsets, contradicting its connectedness. Therefore, both inequalities \eqref{HKS_gamma} and \eqref{HKS_lambda2F_eq} are strict in the class of connected open sets.

It remains to prove sharpness. For all $\varepsilon>0$, let $m_\varepsilon < \frac{V}{2}$ such that $(m_\varepsilon)_{\varepsilon}$ converges to $\frac{V}{2}$, as $\varepsilon \to 0$. For $\varepsilon>0$, we consider $\widetilde{\Omega}_\varepsilon$ a bounded open connected dumbbell containing two disjoint balls of measure $m_\varepsilon$, joined by a tube whose radius tends to $0$, and we choose the construction so that $|\widetilde\Omega_\varepsilon|\to V$. A similar construction appears in \cite[Section 4.1, Proposition 4.1.4]{Henrot}. For $\varepsilon>0$, we set
$$
    s_\varepsilon = \left(\frac{V}{|\widetilde\Omega_\varepsilon|}\right)^{\frac{1}{N}}, \qquad \Omega_\varepsilon = s_\varepsilon \widetilde{\Omega}_\varepsilon.
$$
Then, for all $\varepsilon > 0$, $|\Omega_\varepsilon| = V$, and $s_\varepsilon \to_{\varepsilon\rightarrow 0} 1$. The two rescaled balls that are parts of $\Omega_\varepsilon$ have measure $s_\varepsilon^N m_\varepsilon \to \frac{V}{2}$. For $\varepsilon>0$, let $\varphi_{1,\varepsilon}$ and $\varphi_{2,\varepsilon}$ be their positive first $(2,q)$-eigenfunctions, normalised in $L^q$ and extended by zero to $\Omega_\varepsilon$. These extensions belong to $H_0^1(\Omega_\varepsilon)$ and have disjoint supports. Consequently, by using the pair $(\varphi_{1,\varepsilon},\varphi_{2,\varepsilon})$ in the definition \eqref{Definition_Gamma} of $\Gamma(\Omega_\varepsilon) = \gamma(\Omega_\varepsilon)$, we obtain
$$
    \gamma(\Omega_\varepsilon) = \Gamma(\Omega_\varepsilon) \leq \lambda_1(B_{s_\varepsilon^N m_\varepsilon}).
$$
The scaling property from Lemma \ref{Homogeneity} yields, for all $\varepsilon>0$,
$$
    \lambda_1(B_{s_\varepsilon^N m_\varepsilon}) = \left( s_\varepsilon \left(\frac{2m_\varepsilon}{V}\right)^{\frac{1}{N}} \right)^{-2+N(1-\frac{2}{q})} \lambda_1(B_{V/2})
    \xrightarrow[\varepsilon\rightarrow 0]{} \ \lambda_1(B_{V/2}).
$$
Moreover, inequality \eqref{HKS_gamma} gives
$
    \gamma(\Omega_\varepsilon) \geq \lambda_1(B_{V/2}).
$
Therefore, we obtain
$
    \gamma(\Omega_\varepsilon) \xrightarrow[\varepsilon \rightarrow 0]{} \lambda_1(B_{V/2}).
$
This proves  property \eqref{ConnectedShapeInfimum}.

For $N\geq 2$, property \eqref{ConnectedShapeInfimum} and strictness imply that the infimum is not attained. For $N = 1$, we show that the infimum is actually attained by every interval of length $V$. Take $\Omega = (0,V)$. The positive first eigenfunctions on $(0,\frac{V}{2})$ and $(\frac{V}{2},V)$, normalised in $L^q$ and extended by zero, belong to $H_0^1(\Omega)$ and form an admissible pair for $\Gamma((0,V))$. Thus
$$
    \gamma \left((0,V)\right) = \Gamma \left((0,V)\right)\leq \lambda_1 \left(\left(0,\frac{V}{2}\right)\right).
$$
Combined with inequality \eqref{HKS_gamma}, this yields
$$
    \gamma \left((0,V)\right) = \lambda_1 \left(\left(0,\frac{V}{2}\right)\right), \qquad \lambda_{2}^F \left((0,V)\right) = 2^{1 - \frac{2}{q}} \lambda_1\left(\left(0,\frac{V}{2}\right)\right).
$$
Thus, the infimum is attained by every interval of length $V$.

\end{proof}

The previous proof suggests that the infimum of $\Gamma$ among connected sets is actually attained for a disconnected set, the union of two disjoint balls, that is not part of the considered class. This motivates the generalisation of the definitions of $\Gamma(D)$ and $\lambda_2^F(D)$ for any bounded set $D \subset \R^N$. Thus, we can still define
$$
    \Gamma(D) = \inf_{(v,w)\in\mathcal A_2(D)} \max \{E(v),E(w)\}
$$
and
$$
    \lambda_2^F(D) = 2^{1 - \frac{2}{q}}\Gamma(D),
$$
and adapt the shape optimisation result for disconnected sets.

\begin{rem}
    A priori, we cannot define a value $\gamma(D)$ for disconnected sets, in the same way as we did for connected sets. Indeed, definition \eqref{Definition_gamma} is self-referencing and we had to prove its existence and its variational characterisation by intersecting the first non-trivial Fu\v{c}\'\i k curve with the diagonal.
\end{rem}

\begin{cor}\label{GeneralShapeOptimisationGamma}
    Let $D\subset\R^N$ be a bounded open set with $|D| = V$. Then
    $$
        \Gamma(D) \geq \lambda_1(B_{V/2}), \qquad \lambda_2^F(D) \geq 2^{1 - \frac{2}{q}}\lambda_1(B_{V/2}).
    $$
    Moreover, if $B_1$ and $B_2$ are two disjoint balls of measure $\frac{V}{2}$, then
    $$
        \Gamma(B_1\cup B_2) = \lambda_1(B_{V/2}), \qquad \lambda_2^F (B_1\cup B_2) = 2^{1 - \frac{2}{q}}\lambda_1(B_{V/2}).
    $$
    Also, the infimum in the definition of $\Gamma(D)$ is attained.

    Conversely, if equality holds in one of the inequalities, then, up to zero-measure sets and translations, $D$ is the union of two disjoint balls of equal measure $\frac{V}{2}$.
\end{cor}

\begin{proof}
The lower bound follows exactly as for the connected case, since the proof does not depend on the connectedness of $D$.
If $D=B_1\cup B_2$, the positive first $(2,q)$-eigenfunctions of the two components, normalised in $L^q$ and extended by zero, form an admissible pair whose energies are both equal to $\lambda_1(B_{V/2})$; this gives the reverse inequality, and thus the desired equality.

The infimum in the definition of $\Gamma(D)$ is achieved: a minimising sequence $(v_n,w_n)_n$ is bounded in $H_0^1(D)^2$, the constraints defining $\mathcal A_2(D)$ pass to the limit under strong $L^q(D)$ and a.e.\ convergence, and $E$ is weakly lower semi-continuous.

We finally show the converse, and suppose that $\Gamma(D) = \lambda_1(B_{V/2})$. Let $(v,w) \in \mathcal{A}_2(D)$ be an optimal pair for $\Gamma(D)$, i.e. such that $\max \{E(v),E(w)\} = \Gamma(D) = \lambda_1(B_{V/2})$. As for the connected case, this forces
$$
    |\{v>0\}| = |\{w>0\}| = \frac{V}{2}, \qquad E(v) = E(w) = \lambda_1(B_{V/2}).
$$
Hence, by Lemma \ref{RearrangingThePositivePart}, the two sets $\{v>0\}$ and $\{w>0\}$ are disjoint equal balls $B_1$ and $B_2$ up to sets of measure zero. Since $|D| = V$, we also have $|D \setminus (B_1\cup B_2)| = 0$. Therefore, $D = B_1 \cup B_2$ up to a set of measure zero.
\end{proof}

\begin{rem}
    The infimum of $\lambda_2^F(D)$ among all bounded sets of measure $V$ is therefore achieved for the set $D = B_1 \cup B_2$ (with $B_1$ and $B_2$ two balls of equal measure $\frac{V}{2}$) and is equal to $2^{1 - \frac{2}{q}}\lambda_1(B_{V/2})$. We may notice that this value actually corresponds to the first $(2,q)$-eigenvalue of the set $B_1 \cup B_2$, as seen in \cite[Corollary 2.3]{Brasco_Franzina}, where the $(p,q)$-eigenvalue problem is tackled for disconnected sets: hence
    $$
        \lambda_1(B_1 \cup B_2) = \lambda_2^F(B_1 \cup B_2) = 2^{1 - \frac{2}{q}}\lambda_1(B_{V/2}).
    $$
\end{rem}

\section{The Fu\v{c}\'\i k-selected $(2,q)$-eigenfunctions on a ball}

We recall that by Theorem \ref{FucikSelectedEigenfunction} we have defined for all $\beta > \lambda_1(\Omega)$, a value $\alpha_\beta$ and an associated Fu\v{c}\'\i k-selected $(2,q)$-eigenfunction $u_\beta$, for all connected bounded open sets $\Omega$. We are interested in the symmetric properties satisfied by $u_\beta$. In the following, we restrict our study to the case of a ball $\Omega = B \subset \R^N$ centred at the origin, with $N\geq 2$.

\subsection{Definition and main theorems}

In this subsection we gather some definitions and theorems related to the foliated Schwarz symmetry (we refer to  \cite{Brock_Solynin}, \cite{Girao_Weth}, \cite{Brock_Croce_Guibé_Mercaldo}, and \cite{Brock_Chiacchio_Croce_Mercaldo}).
We will also mainly consider continuous functions $u:B \rightarrow\R$; for measurable functions, all identities will be understood almost everywhere.

\begin{defn}
    Let $\mathcal{H}_0$ be the collection of open half-spaces $H$ in $\R^N$ such that $0\in\partial H$. For any $H\in\mathcal{H}_0$, we denote by $\sigma_H$ the reflection with respect to $\partial H$, that is,
    $$
        \sigma_H u(x) := u(\sigma_H x) \, , \, \forall x\in\R^N.
    $$
The \textbf{two-point rearrangement} (also known as the \textbf{polarisation}) of $u$  with respect to $H$  is defined by
$$
    u_H(x) = \left\{
    \begin{array}{ll}
        \max\{u(x)\,;\,u(\sigma_H x)\} & \textnormal{, if } x\in H \\
        \min\{u(x)\,;\,u(\sigma_H x)\} & \textnormal{, if } x\not\in H\,.
    \end{array}
    \right.
$$
\end{defn}

We may notice that $u_H=u$ if and only if $u(x)\geq u(\sigma_H x)$ for all $x\in H$. Similarly, $u_H = \sigma_H u$ if and only if $u(x)\leq u(\sigma_H x)$ for all $x\in H$.

The two-point rearrangement of a function satisfies the following fundamental properties, as seen in \cite[Lemma 3.1]{Brock_Croce_Guibé_Mercaldo} (see also \cite[Lemma 5.3]{Brock_Solynin} for the case $\Phi\equiv1$).

\begin{ppt}\label{Properties2Point}
    Let $H\in \mathcal{H}_0$.
    \begin{enumerate}
        \item If $A\in C([0,+\infty)\times\R,\R)$, $u:B\rightarrow\R$ is a measurable function, and $x\mapsto A(|x|,u) \in L^1(B)$, then $x\mapsto A(|x|,u_H) \in L^1(B)$, and
        $$
            \int_{B} A(|x|,u)\,dx = \int_{B} A(|x|,u_H)\,dx.
        $$

        \item If $\Phi\in L^\infty(\R)$, and $u\in W^{1,p}(B)$ for some $p\in[1,+\infty)$, then
        $$
            \int_{B} \Phi(u)|\nabla u|^p\,dx = \int_{B} \Phi(u_H)|\nabla u_H|^p\,dx.
        $$
    \end{enumerate}
\end{ppt}

We now define the foliated Schwarz symmetrisation and then see the relation with the above definition.

\begin{defn}
    Let $u:B\rightarrow\R$ be a measurable function. The \textbf{foliated Schwarz symmetrisation} $u^*$ is defined as the (unique) function satisfying the following two properties:
    \begin{enumerate}
        \item There exists a function $w:(0,+\infty) \times [0,\pi] \rightarrow \R$, which is non-increasing in the second variable, such that
        $$
            u^*(x) = w\left(|x|,\arccos\left(\frac{x_1}{|x|}\right)\right) \,\,,\,\, \textnormal{for a.e. } x=(x_1,x_2,...,x_N)\in B\setminus\{0\}.
        $$
        \item For all $a,b\in\R$ such that $a<b$, and all $r\geq 0$, we have
        $$
            \mathcal{H}^{N-1}\{x\in B \, ; \, a<u(x)\leq b, |x|=r\} = \mathcal{H}^{N-1}\{x\in B \, ; \, a<u^*(x)\leq b, |x|=r\}.
        $$
    \end{enumerate}
    If $u^*$ is continuous, its value at $0$ is determined by continuity.
\end{defn}

\begin{defn}
    Denote by $P_N=(1,0,...,0)$ the North pole of the unit sphere $\mathcal{S}^{N-1}\subset \R^N$.
    \begin{enumerate}
    \item
    A function $u$ is said to be \textbf{foliated Schwarz symmetric} with respect to $P_N$ if $u=u^*$.
    \item
    A function $u$ is said to be \textbf{foliated Schwarz symmetric} with respect to a point $P\in\mathcal{S}^{N-1}$ if there exists a rotation $\rho$ around the origin such that $\rho(P_N)=P$ and $u\circ\rho=u^*$.
\end{enumerate}
\end{defn}

In other words, for $0<r<R$ and $c\in\R$, each restricted superlevel set $\{x\in B:|x|=r,\ u(x)\geq c\}$ is a spherical cap centred at $rP$, possibly empty or the whole sphere. In particular, $u$ is axially symmetric with respect to $\R P$.

We conclude this subsection by stating the following fundamental theorem (see \cite[Theorem 3.1]{Brock_Croce_Guibé_Mercaldo} and \cite[Proposition 3.2]{Saldana_Weth}).

\begin{thm}\label{2PointGivesFoliated}
    Let $B=B_R(0)$ and $u\in C(\overline B)$. Suppose that for every $H\in\mathcal{H}_0$, either $u_H = u$ or $u_H = \sigma_H u$. Then $u$ is foliated Schwarz symmetric with respect to some $P\in\mathbb{S}^{N-1}$.
\end{thm}

\subsection{Symmetry of the Fu\v{c}\'\i k-selected $(2,q)$-eigenfunctions}

In order to prove that the Fu\v{c}\'\i k-selected $(2,q)$-eigenfunctions $u_\beta$ associated to the pair $(\alpha_\beta,\beta) \in \Sigma$ is foliated Schwarz symmetric, we apply a similar proof to the one found in \cite[Lemma 4.2 and Proposition 4.1]{Girao_Weth}. We start by a lemma about commuting polarisation and non-decreasing functions.

\begin{Lemma}\label{TwoPointRearrangement And IncreasingFunctions}
    Let $v:B\rightarrow\R$ be a measurable function, and let $T:\R\rightarrow\R$ be a non-decreasing function. Then, for all $H\in\mathcal{H}_0$, we have $(T(v))_H = T(v_H)$.
\end{Lemma}

\begin{proof}
A non-decreasing $T$ commutes with $\max$ and $\min$; the claim follows from the definition of $v_H$.
\end{proof}

The previous general result applies in particular to the non-decreasing functions $x\mapsto x^+ = \max \{x,0\}$ and $x\mapsto -x^- = -\max \{-x,0\} = \min \{x,0\}$. Combining that with the properties \ref{Properties2Point} yield the following corollary.

\begin{cor}\label{Integrals2Point}
    For all functions $v\in H_0^1(B)$ and all $H\in\mathcal{H}_0$, we have
    \begin{enumerate}
        \item $(v^+)_H = (v_H)^+$; thus,
        $$
            \displaystyle \int_B |v^+|^q \,dx = \int_B |(v_H)^+|^q \,dx, \quad  \textnormal{and} \quad \displaystyle \int_B |\nabla (v^+)|^2 \,dx = \int_B |\nabla (v_H)^+|^2 \,dx.
        $$
        \item $(-v^-)_H = -(v_H)^-$; thus,
        $$
            \displaystyle \int_B |v^-|^q \,dx = \int_B |(v_H)^-|^q \,dx, \quad  \textnormal{and} \quad \displaystyle \int_B |\nabla (v^-)|^2 \,dx = \int_B |\nabla (v_H)^-|^2 \,dx.
        $$
    \end{enumerate}
\end{cor}

For all $\beta>\lambda_1(B)$, the Fu\v{c}\'\i k-selected $(2,q)$-eigenfunction $u_\beta$ had been obtained as a limit of a (sub)sequence of minimisers $u_{\beta,\varepsilon}$ of $f_{\beta,\varepsilon}$ on $M^{\beta,\varepsilon}$ as $\varepsilon$ tends to $0$, as done in Propositions \ref{ConvergenceSequence} and \ref{PDE_limit}. In order to show that the limit $u_\beta$ is foliated Schwarz symmetric, we will first show in the following Lemmata that polarising $u_\beta$ is equivalent to passing to the limit the sequence of the polarisations $\left( (u_{\beta,\varepsilon})_H \right)_\varepsilon$.

\begin{Lemma}\label{PolarisationIsAMinimiser}
    Let $\beta > \lambda_1(B)$ and $\varepsilon > 0$. Let $H\in\mathcal{H}_0$, and let $u_{\beta,\varepsilon}$ be a minimiser of $f_{\beta,\varepsilon}$ on $M^{\beta,\varepsilon}$, as obtained in Proposition \ref{InfAchieved}. Then, $(u_{\beta,\varepsilon})_H$ is also a minimiser of $f_{\beta,\varepsilon}$ on $M^{\beta,\varepsilon}$.
\end{Lemma}

\begin{proof}
We refer to the beginning of subsection \ref{Section First Curve: Defining The Family} for a recall of the definitions of the functionals $g_{\beta,\varepsilon}$, $G_{\beta,\varepsilon}$, and $f_{\beta,\varepsilon}$. By using Properties \ref{Properties2Point} and their Corollary \ref{Integrals2Point}, we have
$$
    \norm{((u_{\beta,\varepsilon})_H)^-}_{L^q(B)} = \norm{u_{\beta,\varepsilon}^-}_{L^q(B)} = 1\,,\quad
    \int_B |\nabla ((u_{\beta,\varepsilon})_H)^\pm|^2 \,dx = \int_B |\nabla u_{\beta,\varepsilon}^\pm|^2 \,dx,
$$
and
$$
    \int_B G_{\beta,\varepsilon}((u_{\beta,\varepsilon})_H) \,dx = \int_B G_{\beta,\varepsilon}(u_{\beta,\varepsilon}) \,dx.
$$
Thus, since $((u_{\beta,\varepsilon})_H)^+ = (u_{\beta,\varepsilon}^+)_H$, we have
$$
    \int_B g_{\beta,\varepsilon}((u_{\beta,\varepsilon})_H) ((u_{\beta,\varepsilon})_H)^+ \,dx = \int_B g_{\beta,\varepsilon} (u_{\beta,\varepsilon}) u_{\beta,\varepsilon}^+ \,dx.
$$
Consequently,
$
    f_{\beta,\varepsilon}((u_{\beta,\varepsilon})_H) = f_{\beta,\varepsilon}(u_{\beta,\varepsilon}),
$
and
$
    f'_{\beta,\varepsilon}((u_{\beta,\varepsilon})_H) [((u_{\beta,\varepsilon})_H)^+] = f'_{\beta,\varepsilon}(u_{\beta,\varepsilon}) [u_{\beta,\varepsilon}^+]=0.
$
The second equality implies that $(u_{\beta,\varepsilon})_H \in M^{\beta,\varepsilon}$. Since $f_{\beta,\varepsilon}((u_{\beta,\varepsilon})_H) = f_{\beta,\varepsilon}(u_{\beta,\varepsilon})$ and $u_{\beta,\varepsilon}$ is a minimiser of $f_{\beta,\varepsilon}$ on $M^{\beta,\varepsilon}$, then $(u_{\beta,\varepsilon})_H$ has that same property.
\end{proof}

\begin{Lemma}\label{NormalisedPolarisations}
    Let $\beta > \lambda_1(B)$, and let $\bar{u}_\beta \in H_0^1(B)$ be the strong limit in $H_0^1(B)$ of the sequence of the normalised sequence of minimisers $(\bar{u}_{\beta,\varepsilon})_\varepsilon$ given in Proposition \ref{ConvergenceSequence}, that is,
    $$
        \bar{u}_{\beta,\varepsilon} = - u_{\beta,\varepsilon}^- + \frac{u_{\beta,\varepsilon}^+}{\norm{u_{\beta,\varepsilon}^+}_{L^q(B)}} \xrightarrow[\varepsilon \rightarrow 0]{H_0^1(B)} \bar{u}_\beta.
    $$
    Then, for every $H\in\mathcal{H}_0$, the two-point rearrangement $(\bar{u}_\beta)_H$ is the strong limit in $H_0^1(B)$ of the sequence $((\bar{u}_{\beta,\varepsilon})_H)_\varepsilon$ of the two-point rearrangements of the normalised regularised minimisers.
\end{Lemma}

\begin{proof}
Let $H\in\mathcal{H}_0$, and let $\varepsilon > 0$. We define the following increasing function
$$
    \forall s\in\R, \quad T_\varepsilon(s) = \left\{
    \begin{array}{cl}
        \frac{s}{\norm{u_{\beta,\varepsilon}^+}_{L^q(B)}} & \textnormal{, if } s\geq 0 \\
        s & \textnormal{, if } s < 0\,.
    \end{array}
    \right.
$$
Thus, we have $\bar{u}_{\beta,\varepsilon} = T_\varepsilon (u_{\beta,\varepsilon})$. Hence, by Lemma \ref{TwoPointRearrangement And IncreasingFunctions} and Corollary \ref{Integrals2Point}, we can write
$$
    (\bar{u}_{\beta,\varepsilon})_H = T_\varepsilon ((u_{\beta,\varepsilon})_H) = - ((u_{\beta,\varepsilon})_H)^- + \frac{((u_{\beta,\varepsilon})_H)^+}{\norm{((u_{\beta,\varepsilon})_H)^+}_{L^q(B)}}.
$$
By Lemma \ref{PolarisationIsAMinimiser}, $(u_{\beta,\varepsilon})_H$ is a minimiser of $f_{\beta,\varepsilon}$ on $M^{\beta,\varepsilon}$, and thus $(\bar{u}_{\beta,\varepsilon})_H$ is the nodal normalisation associated to the minimiser $(u_{\beta,\varepsilon})_H$.

Finally, polarisation is $1$-Lipschitz in $L^2(B)$, since
$$
|\max\{a,b\}-\max\{c,d\}|^2+|\min\{a,b\}-\min\{c,d\}|^2\leq|a-c|^2+|b-d|^2
$$
for all real $a,b,c,d$; hence $(\bar u_{\beta,\varepsilon})_H\to(\bar u_\beta)_H$ in $L^2(B)$. Since $(\bar u_{\beta,\varepsilon})_H$ is bounded in $H_0^1(B)$ and $\int_B|\nabla(\bar u_{\beta,\varepsilon})_H|^2\,dx=\int_B|\nabla\bar u_{\beta,\varepsilon}|^2\,dx\to\int_B|\nabla\bar u_\beta|^2\,dx=\int_B|\nabla(\bar u_\beta)_H|^2\,dx$ by Corollary \ref{Integrals2Point}, the convergence holds in $H_0^1(B)$.
\end{proof}

\begin{Lemma}\label{NonnormalisedPolarisationsEquation}
    Let $\beta>\lambda_1(B)$, and let $u_\beta$ be the strong limit in $H_0^1(B)$ of the sequence of minimisers $(u_{\beta,\varepsilon})_{\varepsilon}$, given in Proposition \ref{PDE_limit}. Then, for every $H\in\mathcal{H}_0$, the polarisation $(u_\beta)_H$ is also obtained as the strong limit of minimisers $((u_{\beta,\varepsilon})_H)_\varepsilon$, and satisfies the same Fu\v{c}\'\i k equation:
    $$
        -\Delta (u_\beta)_H = -\alpha_\beta \norm{(u_\beta)_H^-}_{L^q(B)}^{2-q} ((u_\beta)_H^-)^{q-1} + \beta \norm{(u_\beta)_H^+}_{L^q(B)}^{2-q} ((u_\beta)_H^+)^{q-1} \quad \textnormal{in } B.
    $$
\end{Lemma}

\begin{proof}
Let $H\in\mathcal{H}_0$. For all $\varepsilon>0$, we define $\bar{u}_{\beta,\varepsilon} = - u_{\beta,\varepsilon}^- + \frac{u_{\beta,\varepsilon}^+}{\norm{u_{\beta,\varepsilon}^+}_{L^q(B)}}$. By Lemma \ref{NormalisedPolarisations}, the sequence $((\bar{u}_{\beta,\varepsilon})_H)_\varepsilon$ converges strongly in $H_0^1(B)$ to the function $(\bar{u}_\beta)_H$. By Proposition \ref{PDE_limit} and Corollary \ref{Integrals2Point}, we have
$$
    \norm{(u_{\beta,\varepsilon})_H^+}_{L^q(B)} = \norm{u_{\beta,\varepsilon}^+}_{L^q(B)} \quad \xrightarrow[\varepsilon\rightarrow 0]{} \quad \norm{u_\beta^+}_{L^q(B)} = \norm{(u_\beta)_H^+}_{L^q(B)} \, \in \,\, (0,+\infty).
$$
Therefore,
\[
\begin{split}
    (u_{\beta,\varepsilon})_H
    &=-(\bar u_{\beta,\varepsilon})_H^- + \norm{(u_{\beta,\varepsilon})_H^+}_{L^q(B)}(\bar u_{\beta,\varepsilon})_H^+\\
    &\xrightarrow[\varepsilon\to0]{H_0^1(B)} - (\bar u_\beta)_H^-+\norm{(u_\beta)_H^+}_{L^q(B)}(\bar u_\beta)_H^+ = (u_\beta)_H.
\end{split}
\]
Hence, the sequence $((u_{\beta,\varepsilon})_H)_\varepsilon$ converges strongly to $(u_\beta)_H$ in $H_0^1(B)$.

By Lemma \ref{PolarisationIsAMinimiser}, for all $\varepsilon>0$, $(u_{\beta,\varepsilon})_H$ is a minimiser of $f_{\beta,\varepsilon}$ on $M^{\beta,\varepsilon}$. Thus, by Proposition \ref{PDEepsilon},
$$
    \Delta (u_{\beta,\varepsilon})_H + \left[\int_B G_{\beta,\varepsilon}((u_{\beta,\varepsilon})_H)\,dx\right]^{\frac{2}{q}-1} g_{\beta,\varepsilon}((u_{\beta,\varepsilon})_H) - \alpha_{\beta,\varepsilon} ((u_{\beta,\varepsilon})_H^-)^{q-1} = 0 \,\,\,\,\textnormal{in } B,
$$
where $\alpha_{\beta,\varepsilon} = \int_B |\nabla (u_{\beta,\varepsilon})_H^-|^2\,dx = \int_B |\nabla u_{\beta,\varepsilon}^-|^2\,dx$. By letting $\varepsilon$ go to $0$ as in Proposition \ref{PDE_limit}, and since the sequence $((u_{\beta,\varepsilon})_H)_\varepsilon$ converges strongly to $(u_\beta)_H$ in $H_0^1(B)$, we obtain the desired Fu\v{c}\'\i k equation satisfied by $(u_\beta)_H$:
$$
    -\Delta (u_\beta)_H = -\alpha_\beta \norm{(u_\beta)_H^-}_{L^q(B)}^{2-q} ((u_\beta)_H^-)^{q-1} + \beta \norm{(u_\beta)_H^+}_{L^q(B)}^{2-q} ((u_\beta)_H^+)^{q-1},
$$
where $\alpha_{\beta} = \int_B |\nabla (u_{\beta})_H^-|^2\,dx = \int_B |\nabla u_{\beta}^-|^2\,dx$.
\end{proof}

Having obtained the same Fu\v{c}\'\i k equation for the polarisation $(u_\beta)_H$, with the same coefficients, we may now show the main theorem of this section, that the Fu\v{c}\'\i k-selected eigenfunction $u_\beta$ is foliated Schwarz symmetric.

\begin{thm}
    Let $B\subset\R^N$, $N\geq2$, be a ball centred at the origin. For all $\beta > \lambda_1(B)$, the Fu\v{c}\'\i k-selected $(2,q)$-eigenfunction $u_\beta$ is foliated Schwarz symmetric.
\end{thm}

\begin{proof}
By Proposition \ref{PDE_limit}, the $(2,q)$-eigenfunction $u_\beta$ satisfies the following equation
$$
    \left\{
    \begin{array}{cl}
        -\Delta u_\beta = - \alpha_\beta (u_\beta^-)^{q-1} + \beta \norm{u_\beta^+}_{L^q(B)}^{2-q} (u_\beta^+)^{q-1} & \textnormal{ in } B \\
        u_\beta = 0 & \textnormal{ on } \partial B.
    \end{array}
    \right.
$$
Let $H\in\mathcal{H}_0$. By Lemma \ref{NonnormalisedPolarisationsEquation}, $(u_\beta)_H$ satisfies the same equation as $u_\beta$:
$$
    \left\{
    \begin{array}{cl}
        -\Delta (u_\beta)_H = -\alpha_\beta \norm{(u_\beta)_H^-}_{L^q(B)}^{2-q} ((u_\beta)_H^-)^{q-1} + \beta \norm{(u_\beta)_H^+}_{L^q(B)}^{2-q} ((u_\beta)_H^+)^{q-1} & \textnormal{ in } B \\
        (u_\beta)_H = 0 & \textnormal{ on } \partial B,
    \end{array}
    \right.
$$
where $\norm{(u_\beta)_H^+}_{L^q(B)} = \norm{u_\beta^+}_{L^q(B)}$, by Corollary \ref{Integrals2Point}. Moreover, $\sigma_H u_\beta$ also satisfies the same equation as $u_\beta$, because $B$ is a ball centred at the origin and the Laplacian commutes with isometries (such as $\sigma_H$); hence, $\sigma_H u_\beta$ satisfies
$$
    \left\{
    \begin{array}{cl}
        -\Delta (\sigma_H u_\beta) = - \alpha_\beta ((\sigma_H u_\beta)^-)^{q-1} + \beta \norm{u_\beta^+}_{L^q(B)}^{2-q} ((\sigma_H u_\beta)^+)^{q-1} & \textnormal{ in } B \\
        \sigma_H u_\beta = 0 & \textnormal{ on } \partial B.
    \end{array}
    \right.
$$
We may define the following continuous and strictly increasing function
\[
\begin{split}
    \begin{array}{rcl}
        F : \R & \rightarrow & \R \\
        s & \mapsto & - \alpha_\beta (s^-)^{q-1} + \beta \norm{u_\beta^+}_{L^q(B)}^{2-q} (s^+)^{q-1},
    \end{array}
\end{split}
\]
such that all three functions $u_\beta$, $\sigma_H u_\beta$, and $(u_\beta)_H$ satisfy the same equation
$$
    \left\{
    \begin{array}{cl}
        -\Delta v = F(v) & \textnormal{ in } B \\
        {v} = 0 & \textnormal{ on } \partial B.
    \end{array}
    \right.
$$
By Proposition \ref{FucikRegularity}, since the ball $B$ is an open bounded set of class $C^{\infty}$, all functions $u_\beta$, $\sigma_H u_\beta$, and $(u_\beta)_H$ are of class $C^{2,q-1}(\bar{B})$.

Define on $B \cap H$ the function $w = |u_\beta - \sigma_H u_\beta| = 2 (u_\beta)_H - (u_\beta + \sigma_H u_\beta) \geq 0$. The regularity of $u_\beta$, $\sigma_H u_\beta$, and $(u_\beta)_H$ gives us in particular that $w\in C^2(B \cap H) \cap C(\overline{B \cap H})$. Moreover, the function $w$ satisfies the following equation
$$
    \left\{
    \begin{array}{cl}
        -\Delta w = 2 F((u_\beta)_H) - F(u_\beta) - F(\sigma_H u_\beta) & \textnormal{ in } B \cap H \\
        w = 0 & \textnormal{ on } \partial \left(B \cap H\right).
    \end{array}
    \right.
$$
Since the function $F$ is strictly increasing and by the definition of $u_H$ on $B \cap H$, we can write
$$
    -\Delta w = |F(u_\beta) - F(\sigma_H u_\beta)| \geq 0 \quad \textnormal{ in } B \cap H.
$$
Since $w\in C^2(B \cap H) \cap C(\overline{B \cap H})$, the strong maximum principle (see \cite[Section 6.4.2, Theorem 3]{Evans}) implies that either $w=0$ on $B \cap H$ or $w>0$ on $B \cap H$. The former case leads to $u_\beta = \sigma_H u_\beta = (u_\beta)_H$. In the latter, the continuous function $u_\beta - \sigma_H u_\beta$ never vanishes on the connected set $B \cap H$, hence it has constant sign there. This means that either $u_\beta > \sigma_H u_\beta$ on $B \cap H$, which means $(u_\beta)_H = u_\beta$, or that $u_\beta < \sigma_H u_\beta$ on $B \cap H$, which means $(u_\beta)_H = \sigma_H u_\beta$. Since this is true for all $H\in\mathcal{H}_0$, Theorem \ref{2PointGivesFoliated} gives that the eigenfunction $u_\beta$ is foliated Schwarz symmetric with respect to some $P \in \mathbb{S}^{N-1}$.
\end{proof}

In particular, the above applies to the diagonal solution $u_\gamma$.
\begin{cor}\label{ConditionalSecondEigenfunctionSymmetry}
    If the diagonal Fu\v{c}\'\i k-selected $(2,q)$-eigenfunction satisfies $t_\gamma = 1$, then the corresponding $(2,q)$-eigenvalue $\lambda_2^F(B) = 2^{1 - \frac{2}{q}}\gamma$ admits an eigenfunction $u_\gamma$ that is foliated Schwarz symmetric.
\end{cor}

\section*{Acknowledgements}
This work is partially supported by the project STOIQUES ``Shape and Topology
Optimization: Impactful QUestions and Emerging Subjects'' (ANR-24-CE40-2216),
financed by the French Agence Nationale de la Recherche (ANR).
GP gratefully acknowledges the hospitality of SAMM, Universit\'e Paris~1
Panth\'eon-Sorbonne, during his visit, supported by this project.
EA gratefully acknowledges the hospitality of the Dipartimento di Matematica e
Fisica of the Universit\`a degli Studi della Campania ``Luigi Vanvitelli''
during his visit.

GP is a member of the Gruppo Nazionale per l'Analisi Matematica, la
Probabilit\`a e le loro Applicazioni (GNAMPA) of the Istituto Nazionale di Alta
Matematica (INdAM) and wishes to acknowledge financial support from INdAM GNAMPA
Project 2026 ``Problemi di ottimizzazione di forma in contesti anisotropi e
non-locali'' (CUP E53C25002010001). The work of GP was partially supported by
the Portuguese government through FCT -- Funda\c{c}\~ao para a Ci\^encia e a
Tecnologia, I.P., project 2024.14494.PEX with DOI identifier
10.54499/2024.14494.PEX (project ASSO).

\section*{Use of AI tools}
During the preparation of this manuscript the authors used AI tools to proofread the text and to improve its language and readability. After using this tool, the authors reviewed and edited the content as needed and take full responsibility for the content of the publication.

\end{document}